\documentclass[twoside]{amsbook}
\usepackage{fancyheadings}
\usepackage{amssymb}
\usepackage{euscript}
\usepackage{amscd}
\usepackage{mathrsfs}
\usepackage[french]{babel}

\font\titre=cmbx10 scaled 1728
\font\para=cmbx10 scaled 1200

\newcounter{se}
\newcounter{theom}[se]
\gdef\numse{\global\stepcounter{se}%
\these}
\gdef\numtheoreme{\global\stepcounter{theom}%
\thetheom}

\renewcommand{\section}[1]%
{\vskip 7mm plus 20mm minus 1.5mm\penalty-50%
\vskip 0mm plus -20mm minus 1.5mm\penalty-50%
\noindent{\para \numse\ #1}\nobreak}

\renewenvironment{th}[1]%
{\bigskip {{\bf \these.\numtheoreme}\ #1.}~\it}%

\newenvironment{th*}[1]%
{\bigskip {\bf #1.}~\it}%

\renewcommand{\subsection}%
{\bigskip {\bf \these.\numtheoreme} \ }%

\newenvironment{Preuve}{\proof[\it Preuve]}{\qed\endtrivlist}

\newcommand{\id}{\operatorname{id}}

\newcommand{\ad}{\operatorname{{ad}}}

\newcommand{\tr}{\operatorname{{tr}}}
\newcommand{\rg}{\operatorname{{rg}}}

\newcommand{\Mat}{\operatorname{{Mat}}}

\newcommand{\Spec}{\operatorname{{Spec}}}

\renewcommand{\L}{\operatorname{L}}

\newcommand{\Fract}{\operatorname{Fract}}

\renewcommand{\S}{\operatorname{S}}

\newcommand{\U}{\operatorname{U}}

\begin{document}
\setcounter{page}{1}
\thispagestyle{empty}

\cfoot[{}]{{}}
\lhead[\thepage]{{}}
\rhead[{}]{\thepage}

\vspace*{15mm}

\begin{center}
{\titre Alg\`{e}bres de Poisson}\\
\strut\phantom{\small e}\\
{\titre et alg\`{e}bres de Lie r\'{e}solubles}
\end{center}

\medskip

\begin{center}
{\bf Patrice TAUVEL et Rupert W.T. YU}
\end{center}

\vspace{8mm}

{\small {\bf R\'{e}sum\'{e}.}~Soient $\mathfrak{g}$ une alg\`{e}bre de
Lie r\'{e}soluble et $Q$ un id\'{e}al premier $(\ad
\mathfrak{g})$-stable de l'alg\`{e}bre sym\'{e}trique $\S
(\mathfrak{g})$ de $\mathfrak{g}$. Si $E$ est l'ensemble des
\'{e}l\'{e}ments non nuls de $\S (\mathfrak{g})/Q$ qui sont vecteurs
propres pour l'action adjointe de $\mathfrak{g}$ dans $\S
(\mathfrak{g})/ Q$, l'alg\`{e}bre localis\'{e}e $\big( \S
(\mathfrak{g})/ Q\big)_{E}$ a une structure naturelle d'alg\`{e}bre de
Poisson. On \'{e}tudie ici cette alg\`{e}bre.} 

\section{Introduction}

\bigskip

Soient $\Bbbk$ un corps commutatif alg\'{e}briquement clos de
caract\'{e}ristique nulle, $\mathfrak{g}$ une $\Bbbk$-alg\`{e}bre de
Lie r\'{e}soluble, $\U (\mathfrak{g})$ son alg\`{e}bre enveloppante,
et $\S (\mathfrak{g}$ son alg\`{e}bre sym\'{e}trique. On note $\Spec
\U (\mathfrak{g})$ l'ensemble des id\'{e}aux premiers de $\U
(\mathfrak{g})$ et $\Spec \S (\mathfrak{g})^{\mathfrak{g}}$ l'ensemble
des id\'{e}aux premiers $\ad \mathfrak{g}$-stables de $\S
(\mathfrak{g})$. Il existe une bijection canonique $\beta$ de $\Spec
\S (\mathfrak{g})^{\mathfrak{g}}$ sur $\Spec \U (\mathfrak{g})$,
appel\'{e}e la bijection de Dixmier relative \`{a} $\mathfrak{g}$. 

Soient $P \in \Spec \U (\mathfrak{g})$ et $F$ l'ensemble des
\'{e}l\'{e}ments non nuls de $A(P) = \U (\mathfrak{g})/ P$ qui sont
vecteurs propres pour l'action adjointe de $\mathfrak{g}$ dans $A
(P)$. Alors $F$ permet un calcul de fractions dans
$A(P)$. L'alg\`{e}bre associative $A(P)_{F}$ a \'{e}t\'{e}
\'{e}tudi\'{e}e par J.C. Mc\,Connell dans plusieurs articles.

Soient maintenant $Q \in \Spec \S (\mathfrak{g})^{\mathfrak{g}}$ et
$B(Q) = \S (\mathfrak{g})/ Q$. Notons $E$ l'ensemble des vecteurs
propres non nuls pour la repr\'{e}sentation adjointe de $\mathfrak{g}$
dans $B(Q)$. Alors $B(Q)_{E}$ a une structure naturelle d'alg\`{e}bre
de Poisson. Nous \'{e}tudions ici cette alg\`{e}bre de Poisson.

Supposons que $P \in \Spec \U (\mathfrak{g})$ et $Q \in \Spec \S
(\mathfrak{g})^{\mathfrak{g}}$ v\'{e}rifient $P = \beta (Q)$. Les
r\'{e}sultats que l'on obtient pour l'alg\`{e}bre de Poisson
$B(Q)_{E}$ sont plus que voisins de ceux obtenus par
J.C. Mc\,Connell pour l'alg\`{e}bre associative $A(P)_{F}$. D'autre
part, les preuves donn\'{e}es ici suivent de tr\`{e}s pr\`{e}s celles
de Mc\,Connell car, pour l'instant il ne semble exister aucun
r\'{e}sultat reliant la structure de Poisson de $B(Q)_{E}$ et la
structure de l'alg\`{e}bre associative $A(P)_{F}$.

\section{Notations}
\chead[\sl Alg\`{e}bres de Poisson et alg\`{e}bres de Lie
r\'{e}solubles]{\sl 1~~Notations}

\subsection Soit $A$ un anneau commutatif int\`{e}gre. On note $\Fract
A$ le corps des fractions de $A$.

On dit qu'une partie $S$ de $A$ permet un calcul de fractions dans $A$
si elle v\'{e}rifie les deux conditions suivantes :

(i) $1 \in S$ et $0 \notin S$. 

(ii) Le produit de deux \'{e}l\'{e}ments de $S$ appartient \`{a} $S$.

Si elles sont r\'{e}alis\'{e}es, l'ensemble des
\'{e}l\'{e}ments de $\Fract A$ qui s'\'{e}\-cri\-vent $as^{-1}$,
avec $a \in A$ et $s \in S$, est un sous-anneau de $\Fract A$. On le
note $A_{S}$, et on dit que c'est le localis\'{e} de $A$ par $S$. S'il
existe un \'{e}l\'{e}ment non nul $e$ de $A$ tel que $S = \{ e^{n}\, ;
\, n \in \mathbb{N}\}$, on \'{e}crit $A_{e}$ pour $A_{S}$.

\subsection Dans toute la suite, $\Bbbk$ est un corps commutatif de
caract\'{e}ristique nulle. Tous les espaces vectoriels et alg\`{e}bres
consid\'{e}r\'{e}s sont d\'{e}finis sur $\Bbbk$.

Soit $V$ un espace vectoriel. On note $\L (V)$ l'alg\`{e}bre des
endomorphismes de $V$ et $\id_{V}$ l'application identique de $V$. Si
$V$ est de dimension finie, si $\mathcal{B}$ est une base de $V$, et
si $u \in \L (V)$, on d\'{e}signe par $\Mat (u,\mathcal{B})$ la
matrice de $u$ dans la base $\mathcal{B}$.

Si $\mathbb{K}$ est une extension de $\Bbbk$, on d\'{e}signe par $\deg
\tr_{\Bbbk} \mathbb{K}$ le degr\'{e} de transcendance de $\mathbb{K}$
sur $\Bbbk$.

\section{G\'{e}n\'{e}ralit\'{e}s}

\chead[Alg\`{e}bres de Poisson et alg\`{e}bres de Lie r\'{e}solubles]{\sl
2~~G\'{e}n\'{e}ralit\'{e}s} 

\begin{th}{\bf D\'{e}finition}On appelle alg\`{e}bre de Poisson,
ou {\em P}-alg\`{e}bre, une $\Bbbk$-alg\`{e}bre commutative, associative
et unitaire $A$ munie d'une application bilin\'{e}aire altern\'{e}e 
$$
A \times A \to A \ , \ (p,q) \to \{ p, q\} ,
$$
appel\'{e}e crochet de Poisson sur $A$, et v\'{e}rifiant, pour
tous $p,q,r \in A$, les conditions suivantes :

{\em (i)} $ \{ p, \{ q,r\}\} + \{ q,\{ r,p\}\} + \{ r, \{ p,q\} \} =
0$.

{\em (ii)} $\{ pq,r\} = \{ p,r\} q + p \{ q,r\}$.\end{th}

\subsection Soient $A$ et $B$ des P-alg\`{e}bres.

On a $\{ \lambda , p\} = 0$ pour tout $\lambda \in \Bbbk$ et tout $p \in
A$. D'autre part, munie du crochet $\{ .,.\}$, $A$ est une
$\Bbbk$-alg\`{e}bre de Lie.

Une application $f \colon A \to B$ est appel\'{e}e un homomorphisme
de P-alg\`{e}bres si c'est un homomorphisme 
des alg\`{e}bres associatives sous-jacentes \`{a} $A$ et $B$, et si
$$
f(\{ p,q\} ) = \{ f(p), f(q)\}
$$
pour tous $p,q \in A$. On d\'{e}finit de mani\`{e}re \'{e}vidente la
notion d'isomorphisme de P-alg\`{e}bres, ainsi que les notions
d'endomorphisme et d'automorphisme. 

\subsection Jusqu'en 2.8, $A$ est une P-alg\`{e}bre.

Une partie $B$ de $A$ est appel\'{e}e une P-sous-alg\`{e}bre de $A$ si
c'est une sous-alg\`{e}bre associative de $A$ telle que $\{ p,q\} \in
B$ pour tous $p,q \in B$.

Soit $C$ une partie de $A$. Le centralisateur $Y_{A}(C)$ de $C$ dans
$A$ est l'ensemble des $p \in A$ qui v\'{e}rifient $\{ p,q\} = 0$ pour
tout $q \in C$. Il est imm\'{e}diat que $Y_{A}(C)$ est une
P-sous-alg\`{e}bre de $A$. L'ensemble $Y_{A}(A)$, not\'{e} encore
$Y(A)$, est appel\'{e} le centre de $A$. On dit que $A$ est
commutative si $Y_{A}(A) = A$.

\subsection Une partie $J$ de $A$ est appel\'{e}e un P-id\'{e}al de
$A$ si c'est un id\'{e}al de l'alg\`{e}bre associative $A$ et si $\{
p,q\} \in J$ pour tout $(p,q) \in A \times J$. On dit que $A$ est une
P-alg\`{e}bre simple si ses seuls P-id\'{e}aux sont $\{ 0\}$ et $A$.

Soient $J$ un P-id\'{e}al de $A$ et $p \to \overline{p}$ la surjection
canonique $A \to A/J$. Il est imm\'{e}diat de v\'{e}rifier que l'on
munit $A/J$ d'une structure de P-alg\`{e}bre en convenant que
$\{\overline{p}, \overline{q}\} = \overline{\{ p,q\} }$ pour tous $p,q
\in A$. La P-alg\`{e}bre ainsi obtenue est appel\'{e}e la
P-alg\`{e}bre quotient de $A$ par $J$.

\subsection Une P-d\'{e}rivation $\delta$ de $A$ est une
d\'{e}rivation de l'alg\`{e}bre associative $A$ telle que
$$
\delta (\{ p,q\} ) = \{ \delta (p),q\} + \{ p, \delta (q)\}
$$
pour tous $p,q \in A$. 

Si $a \in A$, l'application
$$
d_{a} \colon A \to A \ , \ p \to \{ a,p\}
$$
est une P-d\'{e}rivation de $A$. On dit que c'est la P-d\'{e}rivation
int\'{e}rieure de $A$ d\'{e}finie par $a$.

\subsection Soient $\delta$ une P-d\'{e}rivation de $A$ et
$X$ une ind\'{e}termin\'{e}e sur $A$. Il existe une et une seule
structure de P-alg\`{e}bre sur l'alg\`{e}bre associative $A[X]$, qui
prolonge la structure de $P$-alg\`{e}bre de $A$, et telle que
$$
\{ X, p\} = \delta (p)
$$ 
pour tout $p \in A$. La P-alg\`{e}bre ainsi obtenue est not\'{e}e
$A_{\delta}\{ X\}$. Si $\delta = 0$, on \'{e}crit $A\{ X\}$ pour
$A_{\delta}\{ X\}$.

Plus g\'{e}n\'{e}ralement, si $B$ est une P-sous-alg\`{e}bre de $A$,
et si $x \in Y_{A}(B)$, on d\'{e}signe par $B\{ x\}$ la
P-sous-alg\`{e}bre de $A$ engendr\'{e}e par $B$ et $x$.

\subsection Soient $S$ une partie de $A$
qui permet un calcul des fractions dans $A$ et $A_{S}$ l'alg\`{e}bre
associative localis\'{e}e de $A$ par $S$. Sur $A_{S}$, il existe une
et une seule structure de P-alg\`{e}bre prolongeant celle de $A$. Si
$p,q \in A$ et $s,t \in S$, on a :
$$
\{ ps^{-1}, qt^{-1}\} = \{ p,q\} s^{-1}t^{-1} - \{ p,t\} qs^{-1}t^{-2}
- \{ q,s\} ps^{-2}t^{-1} + \{ s,t\} pqs^{-2}t^{-2}.
$$

\subsection Soient $A$ et $B$ des P-alg\`{e}bres. Il existe une et une
seule structure de P-alg\`{e}bre sur $A \otimes_{\Bbbk} B$, prolongeant
celles de $A$ et $B$, et telle que $\{ p, q\} = 0$ pour tout $(p,q)
\in A \times B$.

\begin{th}{\bf Proposition}Soient $L$ une {\em P}-alg\`{e}bre qui est
un corps, $\delta$ une {\em P}-d\'{e}\-ri\-va\-tion non int\'{e}rieure
de $L$, et $F = L_{\delta} \{ X\}$.

{\em (i)} La {\em P}-alg\`{e}bre $F$ est simple.

{\em (ii)} Le centre de la {\em P}-alg\`{e}bre $\Fract F$ est
l'ensemble des $y \in Y(L)$ qui v\'{e}rifient $\delta (y) =
0$.\end{th}

\begin{Preuve}Si $a \in F$, on note $\deg a$ son degr\'{e} en $X$.

(i) Supposons que $F$ ne soit pas simple, et soit $I$ un P-id\'{e}al
de $F$, non nul et distinct de $F$. Soit $n = \min \{ \deg a \,~; \, a
\in I \backslash \{ 0\} \}$. On a $I \cap \Bbbk = \{ 0\}$ puisque $I \ne
F$, donc $n > 0$. Comme $L$ est un corps, il existe un \'{e}l\'{e}ment
$a$ de $I$ de la forme
$$
X^{n} + a_{n-1} X^{n-1} + \cdots + a_{1} X + a_{0},
$$
avec $a_{0}, \dots , a_{n-1} \in L$. Si $u \in L$, on a :
$$
\{ a,u\} = n X^{n-1} \delta (u) + X^{n-1} \{ a_{n-1}, u\} + b,
$$ 
o\`{u} $b \in F$ v\'{e}rifie $\deg b < n-1$. Comme $\{ a,u\} \in
I$, il r\'{e}sulte du choix de $a$ que $\{ a,u\} = 0$, donc $n \delta
(u) + \{a_{n-1}, u\} = 0$. Par suite, $\delta = -\dfrac{1}{n}
d_{a_{n-1}}$, ce qui est absurde.

(ii) a) Soit $a = X^{m}a_{m} + X^{m-1}a_{m-1} + \cdots + a_{1}X +
a_{0}$ un \'{e}l\'{e}ment central de $F$. Prouvons que $a \in Y(L)$ et
que $\delta (a) = 0$. On a :
$$
0 = \{ X, a\} = \operatornamewithlimits{\textstyle \sum}_{i=0}^{m}
X^{i} \delta (a_{i}).
$$
De m\^{e}me, pour $u \in L$ :
$$
0 = \{ u, a\} = \operatornamewithlimits{\textstyle \sum}_{i=0}^{m}
X^{i}  \{ u, a_{i}\} - \operatornamewithlimits{\textstyle
\sum}_{i=1}^{m} iX^{i-1} \delta (u) a_{i}.
$$
On en d\'{e}duit que $a_{m} \in Y(L)$ et que $\delta (a_{m}) =
0$. Quitte \`{a} remplacer $a$ par $aa_{m}^{-1}$, on peut donc
supposer que $a_{m} = 1$. Les relations pr\'{e}c\'{e}dentes montrent
alors que $m \delta (u) = \{ u,a_{m-1}\}$ pour tout $u \in L$. Comme
$\delta$ n'est pas int\'{e}rieure, on a ainsi $m = 0$.

b) On va prouver qu'un \'{e}l\'{e}ment central de $\Fract F$ est
quotient de deux \'{e}l\'{e}ments centraux de $F$. On aura donc le
r\'{e}sultat d'apr\`{e}s ce qui pr\'{e}c\`{e}de.

Soit $a = bc^{-1}\in Y(\Fract F)$, avec $b,c \in F$. On raisonne par
r\'{e}\-cur\-ren\-ce sur $\alpha (a) = \deg b + \deg c$, le cas o\`{u}
$\alpha (a) = 0$ \'{e}tant clair. Quitte \`{a} changer $a$ en
$a^{-1}$, on peut supposer que $\deg b \geqslant \deg c$. Ecrivons $b
= rc + s$, avec $r,s \in F$ et $\deg s < \deg c$. Si $u \in F$, il
vient :
$$
0 = \{ u, a\} = \{ u,r\} + \{ u , sc^{-1}\} = \{ u, r\} + \{ u, s\} c^{-1}
- s \{ u, c\} c^{-2}.
$$
Par suite :
\begin{equation}\label{1}
c^{2} \{ u, r\} + c\{ u, s\} - s \{u,c\} = 0.
\end{equation}
Supposons $u \in L$. Si $\{ u, r\} \ne 0$, on obtient une
contradiction d'apr\`{e}s \eqref{1} car, 
$$
\deg \big( c^{2} \{u, r\} ) > \deg \big (c \{ u,s\} - s \{ u , c\}
\big) .
$$
De m\^{e}me, si $u = X$ on voit, pour la m\^{e}me raison, que $\{
X,r\} = 0$.

D'apr\`{e}s ce qui pr\'{e}c\`{e}de, on a $r \in Y(F)$, et $e =
sc^{-1}$ appartient au centre de $\Fract F$. Comme $\alpha (e) =
\deg s + \deg c < \alpha (a)$, l'hypoth\`{e}se de r\'{e}currence
montre que $e = s'c'^{-1}$, avec $s', c' \in Y (F)$. Comme $a = (rc' +
s')c'^{-1}$, et que $rc' + s', c'$ sont des \'{e}l\'{e}ments de
$Y(F)$, on a obtenu le r\'{e}sultat. \end{Preuve}

\section{Un exemple}
\chead[Alg\`{e}bres de Poisson et alg\`{e}bres de Lie
r\'{e}solubles]{\sl 3~~Un exemple} 

\subsection Soient $n$ un entier positif ou nul et $X_{1},Y_{1}, \dots ,
X_{n}, Y_{n}$ des ind\'{e}termin\'{e}es sur $\Bbbk$. Sur $A = \Bbbk
[X_{1}, \dots , X_{n}, Y_{1}, \dots , Y_{n}]$, il existe une et une
seule structure de P-alg\`{e}bre telle que, pour $1 \leqslant i,j
\leqslant n$, on ait
$$
\{ X_{i}, X_{j}\} = \{ Y_{i}, Y_{j}\} = 0 \ , \ \{ X_{i}, Y_{j}\} =
\delta_{ij}, 
$$
o\`{u} $\delta_{ij}$ est le symbole de Kronecker. La P-alg\`{e}bre
ainsi obtenue est not\'{e}e $B_{n}$. 

Dans la suite, si $\alpha = (i_{1}, \dots , i_{n})$ et $\beta =
(j_{1}, \dots , j_{n})$ sont des \'{e}l\'{e}ments de $\mathbb{N}^{n}$,
on note $X^{\alpha}Y^{\beta}$ l'\'{e}l\'{e}ment de $A$ d\'{e}fini par :
$$
X^{\alpha}Y^{\beta} = X_{1}^{i_{1}}\cdots X_{n}^{_{i_{n}}}
Y_{1}^{j_{1}} \cdots Y_{n}^{i_{n}}.
$$

\begin{th}{\bf Remarque}{\em Soient $V$ un $\Bbbk$-espace vectoriel de
dimension paire $2n$ et $\omega$ une forme bilin\'{e}aire altern\'{e}e
et non d\'{e}g\'{e}n\'{e}r\'{e}e sur $V$. Il existe une et une seule
structure de P-alg\`{e}bre sur l'alg\`{e}bre sym\'{e}trique $\S (V)$
de $V$ telle que
$$
\{ v, w\} = \omega (v,w)
$$
pour tous $v,w \in V$. La P-alg\`{e}bre ainsi obtenue est isomorphe
\`{a} $B_{n}$}. \end{th}

\subsection Soit $Z$ une $\Bbbk$-alg\`{e}bre associative et commutative,
consid\'{e}r\'{e}e comme une P-alg\`{e}bre commutative. On d\'{e}signe
par $B_{n}(Z)$ la P-alg\`{e}bre $Z \otimes_{\Bbbk} B_{n}$ (voir
2.8). Tout \'{e}l\'{e}ment $a$ de $B_{n}(Z)$ s'\'{e}crit uniquement
sous la forme
$$
a = \operatornamewithlimits{\textstyle \sum}_{\alpha , \beta \in
\mathbb{N}^{n}} \lambda_{\alpha , \beta} X^{\alpha} Y^{\beta},
$$
o\`{u} les $\lambda_{\alpha , \beta}$ sont des \'{e}l\'{e}ments de
$Z$. Si $1 \leqslant i \leqslant n$, on a :
\begin{equation} \label{2}
\{ X_{i},a\} = \frac{\partial a}{\partial Y_{i}} \ \raisebox{2pt}{,} \
\{ Y_{i}, a\} = -\frac{\partial a}{\partial X_{i}} \cdotp
\end{equation}

\begin{th}{\bf Lemme}Soit $B = B_{n}(Z)$ comme pr\'{e}c\'{e}demment.

{\em (i)} Le centre de $B$ est \'{e}gal \`{a} $Z$.

{\em (ii)} Tout {\em P}-id\'{e}al de $B$ est engendr\'{e} par son
intersection avec $Z$.

{\em (iii)} Plus g\'{e}n\'{e}ralement, si $C$ est une {\em
P}-alg\`{e}bre, tout {\em P}-id\'{e}al de $C \otimes_{\Bbbk} B_{n}$ est de
la forme $J \otimes_{\Bbbk} B_{n}$, o\`{u} $J$ est un {\em P}-id\'{e}al de
$C$.

{\em (iv)} Soit $\delta$ une {\em P}-d\'{e}rivation de $B$. Il existe $b
\in B$ tel que $\delta' = \delta - d_{b}$ v\'{e}rifie :
$$ 
\delta'|1 \otimes B_{n} = 0 \ , \ \delta' | Z \otimes 1 = \delta |
Z \otimes 1.
$$\end{th}

\begin{Preuve}(i) C'est imm\'{e}diat en utilisant les relations \eqref{2}.

(iii) Si $\alpha = (i_{1}, \dots , i_{n}) \in \mathbb{N}^{n}$, posons
$|\alpha | = i_{1} + \cdots + i_{n}$, $\alpha ! = i_{1}! \cdots
i_{n}!$, et $(d_{X})^{\alpha} = (d_{X_{1}})^{i_{1}}{\circ} \cdots
{\circ} (d_{X_{n}})^{i_{n}}$.

Soit $a = \lambda \otimes X^{\mu} Y^{\nu}\in B$. Si $|\mu | + |\nu |
\leqslant |\alpha | + |\beta |$, il vient :
$$
\begin{cases}
(d_{X})^{\alpha }{\circ} (d_{Y})^{\beta }(a) = 0 \ \text{ si } (\mu
, \nu ) \ne (\beta , \alpha ) \\ 
(d_{X})^{\alpha} {\circ} (d_{Y})^{\beta}(a) = (-1)^{|\beta |}\alpha !
\beta ! \lambda \ \text{ si } \ (\mu , \nu ) = (\beta , \alpha ).
\end{cases}
$$
La caract\'{e}ristique de $\Bbbk$ \'{e}tant nulle, on en d\'{e}duit
imm\'{e}diatement (iii). L'assertion (ii) en est un cas particulier.

(iv) Comme $Y(B) = Z$, on a $\delta (Z) \subset Z$. Pour $a \in B_{n}$
et $z \in Z$, posons $\delta'(z \otimes a) = \delta (z) \otimes a$. On
v\'{e}rifie facilement que $\delta'$ est une P-d\'{e}rivation de
$B$. On est donc ramen\'{e} \`{a} prouver que, si $\delta$ est une
P-d\'{e}rivation de $B$ telle que $\delta |Z = 0$, alors $\delta$ est
int\'{e}rieure.

Si $1 \leqslant i \leqslant n$, posons $\delta (X_{i}) = q_{i}$ et
$\delta (Y_{i}) = -p_{i}$. D'apr\`{e}s les relations de 3.1, il vient
\begin{align*}
\{ \delta (X_{i}), X_{j}\} + \{ X_{i}, \delta (X_{j})\} & = \{
\delta (Y_{i}), Y_{j}\} + 
\{ Y_{i}, \delta (Y_{j})\} = 0 \\
{} & = \{ \delta (X_{i}), Y_{j}\} + \{ X_{i}, \delta (Y_{j})\} ,
\end{align*}
soit :
$$
- \frac{\partial q_{i}}{\partial Y_{j}} + \frac{\partial
q_{j}}{\partial Y_{i}} = - \frac{\partial p_{i}}{\partial X_{j}} +
\frac{\partial p_{j}}{\partial X_{i}} = \frac{\partial q_{i}}{\partial
X_{j}} - \frac{\partial p_{j}}{\partial Y_{i}} = 0.
$$
Par suite, il existe $b \in B$ tel que, pour $1 \leqslant i \leqslant
n$, on ait :
$$
\frac{\partial b}{\partial X_{i}} = p_{i} \ , \ \frac{\partial
b}{\partial Y_{i}} = q_{i}.
$$
Si $\delta' = d_{-b}$, il vient alors :
$$
\begin{cases}
\delta (X_{i}) = q_{i} = \dfrac{\partial b}{\partial Y_{i}} = \{
X_{i}, b\} = \delta'(X_{i}) , \\
\delta (Y_{i}) = -p_{i} = - \dfrac{\partial b}{\partial X_{i}} = \{
Y_{i}, b\} = \delta'(Y_{i}).
\end{cases}
$$
D'o\`{u} $\delta = \delta'$. \end{Preuve}

\begin{th}{\bf Lemme}{\em (i)} Soient $A$ une {\em P}-alg\`{e}bre, $B$
et $C$ deux {\em P}-sous-alg\`{e}bres permutables de $A$, engendrant
$A$. On suppose qu'il existe un entier $n$ tel que $C$ soit isomorphe
\`{a} $B_{n}$. Alors l'homomorphisme canonique de $B \otimes_{\Bbbk}
C$ dans $A$ est un isomorphisme de {\em P}-alg\`{e}bres.

{\em (ii)} Soient $X,Y$ des ind\'{e}termin\'{e}es sur $A$ et $\Delta$
la {\em P}-d\'{e}rivation de $A\{ Y\}$ telle que $\Delta (Y) = 1$,
$\Delta | A = 0$. La {\em P}-sous-alg\`{e}bre de $(A\{ Y\} )_{\Delta}
\{ X\}$ engendr\'{e}e par $X$ et $Y$ est isomorphe \`{a}
$B_{1}$. L'homomorphisme canonique de $A \otimes_{\Bbbk} B_{1}$ dans $(A
\{ Y\} )_{\Delta}\{ X\}$ est un isomorphisme de {\em
P}-alg\`{e}bres. \end{th}

\begin{Preuve}(i) Soit $\varphi \colon B \otimes_{\Bbbk} C \to A$
l'homomorphisme canonique. Comme $\{ B, C\} = \{ 0\}$, c'est un
homomorphisme de P-alg\`{e}bres ; son noyau est donc un P-id\'{e}al.

Compte tenu de 3.4, tout P-id\'{e}al de $B\otimes_{\Bbbk}C$ est de la
forme $J \otimes_{\Bbbk} C$, o\`{u} $J$ est un P-id\'{e}al de $B$. D'autre
part, la restriction de $\varphi$ \`{a} $B$ est injective. On a donc
$\ker \varphi = \{ 0\}$. L'application $\varphi$ \'{e}tant surjective,
on voit donc que c'est un isomorphisme de P-alg\`{e}bres.

(ii) Le premier point est clair. Le second est un cas particulier de
(i). \end{Preuve}

\begin{th}{\bf Lemme}Soient $A$ une {\em P}-alg\`{e}bre, $\delta$ une {\em
P}-d\'{e}rivation localement nilpotente de $A$, et $\alpha$ un
\'{e}l\'{e}ment central de $A$ tel que $\delta (\alpha ) = 1$. On note $B$
la {\em P}-alg\`{e}bre quotient $A/A \alpha$ et $p \to \overline{p}$
la surjection canonique $A \to B$. Soient $Y$ une ind\'{e}termin\'{e}e
et $\Delta$ la {\em P}-d\'{e}rivation de $B \{ Y\}$ telle que $\Delta
(Y) = 1$, $\Delta | B = 0$. On d\'{e}finit $\chi \colon A \to B \{
Y\}$ en posant, pour $p \in A$ :
$$
\chi (p) = \operatornamewithlimits{\textstyle \sum}_{n \geqslant 0}
\frac{1}{n!} \overline{\delta^{n}(p)} Y^{n}.
$$

Alors $\chi$ est un isomorphisme de {\em P}-alg\`{e}bres tel que $\delta = 
\chi^{-1} {\circ} \Delta {\circ} \chi$. On a $\chi^{-1}(Y) = \alpha$
et, si $p \in A$ :
$$
\chi^{-1}(\overline{p}) = \operatornamewithlimits{\textstyle \sum}_{m
\geqslant 0} \frac{(-1)^{m}}{m!} \delta^{m}(p) \alpha^{m}.
$$ \end{th}

\begin{Preuve}On v\'{e}rifie facilement que $\chi {\circ} \delta = \Delta
{\circ} \chi$ et que $\chi (pq) = \chi (p) \chi (q)$ pour $p,q \in
A$. D'autre part :
$$
\{ \chi (p), \chi (q)\} = \operatornamewithlimits{\textstyle
\sum}_{m,n\geqslant 0} \frac{1}{m!}\frac{1}{n!}
\overline{\{ \delta^{m}(p), \delta^{n}(q)\}} Y^{m+n}.
$$
En remarquant que :
$$
\delta^{m+n}(\{ p,q\} ) = \operatornamewithlimits{\textstyle
\sum}_{k=0}^{m+n} \begin{pmatrix}
n+m \\
k\end{pmatrix} \{ \delta^{k}(p), \delta^{m+n-k}(q)\} ,
$$
on voit que $\{ \chi (p), \chi (q)\} = \chi (\{ p,q\} )$. Ainsi,
$\chi$ est un homomorphisme de P-alg\`{e}bres.

Soit $r \in A$. Il vient :
\begin{align*}
\operatornamewithlimits{\textstyle \sum}_{m \geqslant 0}
\frac{(-1)^{m}}{m!} \delta^{m}(r\alpha ) \alpha^{m} & =
\operatornamewithlimits{\textstyle \sum}_{m \geqslant 0}
\frac{(-1)^{m}}{m!} \big( \delta^{m}(r)\alpha + m\delta^{m-1}(r)\big)
\alpha^{m} \\
{} & = \operatornamewithlimits{\textstyle \sum}_{m \geqslant 0}
\frac{(-1)^{m}}{m!} \delta^{m}(r) \alpha^{m+1} \\
{} & \phantom{abcdef} + 
\operatornamewithlimits{\textstyle \sum}_{m \geqslant 0}
\frac{(-1)^{m+1}}{m!} \delta^{m}(r) \alpha^{m+1} = 0
\end{align*}
On en d\'{e}duit qu'il existe une unique application lin\'{e}aire
$\theta \colon B \{ Y\} \to A$ telle que, pour $p \in A$ et $n \in
\mathbb{N}$ :
$$
\theta (\overline{p}Y^{n}) = \operatornamewithlimits{\textstyle
\sum}_{m \geqslant 0} \frac{(-1)^{m}}{m!} \delta^{n}(p) \alpha^{m+n}.
$$
On v\'{e}rifie facilement que $\theta$ est un homomorphisme de
P-alg\`{e}bres. Il vient :
$$
\theta {\circ} \chi (p) = \operatornamewithlimits{\textstyle
\sum}_{m,n \geqslant 0} \frac{1}{m!} \frac{(-1)^{n}}{n!} \delta^{m+n}(p)
\alpha^{m+n} = p \ , \ \chi {\circ} \theta (Y) = \chi (\alpha ) = Y.
$$
De m\^{e}me :
\begin{align*}
\chi {\circ} \theta (p) & = \chi \Big(
\operatornamewithlimits{\textstyle \sum}_{m\geqslant 0}
\frac{(-1)^{m}}{m!} \delta^{m}(p) \alpha^{m}\Big) =
\operatornamewithlimits{\textstyle \sum}_{m,n \geqslant 0}
\frac{(-1)^{m}}{m!n!} \overline{\delta^{n}\big(
\delta^{m}(p)\alpha^{m}\big)} Y^{n} \\
{} & = \operatornamewithlimits{\textstyle \sum}_{m,k,\ell \geqslant 0}
\frac{(-1)^{m}}{m!k!\ell
!}\overline{\delta^{k+m}(p)\delta^{\ell}(\alpha^{m})} Y^{k+\ell} \\
{} & = \operatornamewithlimits{\textstyle \sum}_{k,m \geqslant 0}
\frac{(-1)^{m}}{m!k!m!} m! \overline{\delta^{k+m}(p)} Y^{k+m} =
\overline{p}. 
\end{align*}
On a donc obtenu le r\'{e}sultat. \end{Preuve}

\begin{th}{\bf Remarque}{\em Avec les hypoth\`{e}ses et notations de
3.6, il est imm\'{e}diat que $\ker \delta$ est une {\em
P}-sous-alg\`{e}bre de $A$ qui est isomorphe \`{a} $B$.} \end{th}

\begin{th}{\bf Lemme}Soit $C$ une {\em P}-alg\`{e}bre, $A$ une {\em
P}-sous-alg\`{e}bre commutative de $C$, et $L$ un sous-espace de
dimension finie de $C$. On suppose v\'{e}rifi\'{e}es les conditions
suivantes :

{\em a)} Pour tout \'{e}l\'{e}ment $u$ de $L$, $d_{u}$ induit une {\em
P}-d\'{e}rivation localement nilpotente de $A$ et, si $x,y \in L$,
on a $d_{x} {\circ} d_{y}| A = d_{y}{\circ} d_{x}|A$.

{\em b)} L'ensemble des $p \in A$ v\'{e}rifiant $d_{x}(p) = 0$ pour
tout $x \in L$ est \'{e}gal \`{a} $\Bbbk$.

Alors il existe un sous-espace de dimension finie $V$ de $A$ tel que :

{\em (i)} $A$ est isomorphe \`{a} l'alg\`{e}bre sym\'{e}trique $\S
(V)$ de $V$ (consid\'{e}r\'{e}e comme {\em P}-alg\`{e}bre
commutative), et les seuls id\'{e}aux de $\S (V)$ stables par les
$d_{x}$, avec $x \in L$, sont $\{ 0\}$ et $\S (V)$.

{\em (ii)} Pour tout $x \in L$ et tout $v \in V$, on a $d_{x}(v) \in
\Bbbk$. \end{th}

\begin{Preuve}Si $n \in \mathbb{N}^{*}$, notons $A_{n}$ l'ensemble
des \'{e}l\'{e}ments $p$ de $A$ tels que $d_{x_{1}} {\circ} \cdots
{\circ} d_{x_{n}}(p) = 0$ pour tous $x_{1}, \dots , x_{n} \in L$. On a
ainsi $A_{1} = \Bbbk$. Pour obtenir le r\'{e}sultat, on peut supposer $A
\ne \Bbbk$.

Soit $p_{1} \in A_{2} \backslash \{ 0\}$. On a $d_{x}(p_{1}) \in \Bbbk$
pour tout $x \in L$, et il existe $x_{1} \in L$ tel que
$d_{x_{1}}(p_{1}) = 1$. Posons $B = \ker d_{x_{1}}$. D'apr\`{e}s 3.6
et 3.7, les P-alg\`{e}bres $A$ et $B\{ X\}$ sont isomorphes. Si $\dim
L = 1$, on a $B = \Bbbk$, et on a obtenu le r\'{e}sultat. Supposons $\dim
L \geqslant 2$, et raisonnons par r\'{e}currence sur la dimension de
$L$.

Soient $x \in L$ et $y = x - d_{x}(p_{1}) x_{1}$. Il vient :
$$
d_{x}(p_{1}) \in \Bbbk \ , \ d_{y}(p_{1}) = 0.
$$
D'autre part, $d_{x}|A$ et $d_{y}|A$ commutent puisque $d_{x}|A$ et
$d_{x_{1}}|A$ commutent. Il en r\'{e}sulte que $B = \ker d_{x_{1}}$
est stable par les $d_{y}$, avec $y \in L'$, o\`{u}
$$
L' = \{x - d_{x}(p_{1}) x_{1} \ ; \, x \in L\}.
$$ 
Enfin, si $y \in L'$, $d_{y}$ induit une d\'{e}rivation localement
nilpotente de $B$, et on a $\dim L' \leqslant \dim L -1$.

D'apr\`{e}s l'hypoth\`{e}se de r\'{e}currence, il existe un
sous-espace de dimension finie $W$ de $B$ v\'{e}rifiant les conditions
suivantes :

(1) En tant que P-alg\`{e}bres commutatives, $B$ et $\S (W)$ sont
isomorphes. 

(2) Pour $y \in L'$ et $w \in W$, on a $d_{y}(w) \in \Bbbk$.

(3) Les seuls id\'{e}aux de $\S (W)$ stables par les $d_{y}$, avec $y
\in L'$, sont $\{ 0\}$ et $\S (W)$.

Or, si $y = x - d_{x}(p_{1})x_{1}$, avec $x \in L$, on a $d_{y}(w) =
d_{x}(w)$ pour $y \in W$. Par suite, $d_{x}(w) \in \Bbbk$. En prenant
$V = W + \Bbbk p_{1}$, on obtient alors facilement le
r\'{e}sultat. \end{Preuve}

\subsection Rappelons que $B_{1}$ est la P-alg\`{e}bre $\Bbbk[X_{1},
Y_{1}]$, avec $\{ X_{1}, Y_{1}\} = 1$.

\begin{th*}{\bf Lemme}On conserve les hypoth\`{e}ses et notations de
{\em 3.6}. Il existe un unique homomorphisme de {\em P}-alg\`{e}bres
$\chi \colon A_{\delta}\{ X \} \to B \otimes_{\Bbbk} B_{1}$ tel que
$$
\chi (X) = 1 \otimes X_{1} \ , \ \chi (p) =
\operatornamewithlimits{\textstyle \sum}_{n \geqslant 0}
\frac{1}{n!}\overline{\delta^{n}(p)} \otimes Y_{1}^{n}
$$
pour tout $p \in A$. Cet homomorphisme est un isomorphisme, et on a 
\begin{gather*}
\chi^{-1}(1 \otimes X_{1}) = X \ , \ \chi^{-1}(1 \otimes Y_{1}) =
\alpha ,\\
\chi^{-1}(\overline{p} \otimes 1) = \operatornamewithlimits{\textstyle
\sum}_{m \geqslant 0} \frac{(-1)^{m}}{m!} \delta^{m}(p) \alpha^{m}
\end{gather*}
pour tout $p \in A$. \end{th*}

\begin{Preuve}D'apr\`{e}s 3.5 et 3.6, on a des isomorphismes de
P-alg\`{e}bres :
$$
A_{\delta}\{ X\} \overset{\chi_{1}}{\longrightarrow} (B\{ Y\} )_{\Delta} \{
X\} \overset{\chi_{2}}{\longrightarrow} B \otimes_{\Bbbk} B_{1}.
$$
Soit $\chi = \chi_{2} {\circ} \chi_{1}$. Alors
\begin{gather*}
\chi (X) = \chi_{2}(X) = 1 \otimes X_{1} \ , \ \chi (\alpha ) =
\chi_{2}(Y) = 1 \otimes Y_{1},\\
\chi (p) = \chi_{2} \Big( \operatornamewithlimits{\textstyle \sum}_{n
\geqslant 0} \frac{1}{n!} \overline{\delta^{n}(p)}Y^{n} \Big) =
\operatornamewithlimits{\textstyle \sum}_{n \geqslant 0}
\overline{\delta^{n}(p)} \otimes Y_{1}^{n}\\
\chi^{-1}(\overline{p} \otimes 1) = \chi_{1}^{-1}(\overline{p}) =
\operatornamewithlimits{\textstyle \sum}_{n \geqslant 0}
\frac{(-1)^{n}}{n!} \delta^{n}(p) \alpha^{n}.
\end{gather*}
D'o\`{u} le lemme. \end{Preuve}

\subsection Avec les notations de 3.1, soient $X = X_{1}\cdots X_{n}$
et $S = \{ X^{p}\, ; p \in \mathbb{N}\}$. Alors $S$ est une partie de
$\Bbbk [X_{1}, Y_{1}, \dots , X_{n}, Y_{n}]$ qui permet un calcul de
fractions. Dans la suite, on note $B'_{n}$ la P-alg\`{e}bre
$(B_{n})_{S}$ d\'{e}finie comme en 2.7. Si $Z$ est une alg\`{e}bre
commmutative et associative, consid\'{e}r\'{e}e comme une
P-alg\`{e}bre commutative, on d\'{e}signe par $B'_{n}(Z)$ la
P-alg\`{e}bre $Z \otimes_{\Bbbk} B'_{n}$.

De m\^{e}me, la P-alg\`{e}bre $\Fract B_{n}$ est not\'{e}e $F_{n}$ et,
avec les notations pr\'{e}\-c\'{e}\-den\-tes, on pose $F_{n}(Z) = Z
\otimes_{\Bbbk} F_{n}$.

\section{Alg\`{e}bres de Lie et alg\`{e}bres de Poisson}
\chead[\sl Alg\`{e}bres de Poisson et alg\`{e}bres de Lie
  r\'{e}solubles]{\sl 4~~Alg\`{e}bres de Lie et alg\`{e}bres de Poisson}

\subsection En ce qui concerne les rappels du paragraphe 4, le lecteur
pourra se reporter \`{a} [2] et [11].

Soient $\mathfrak{h}$ une $\Bbbk$-alg\`{e}bre de Lie de dimension
finie, $\U (\mathfrak{h})$ son alg\`{e}bre enveloppante, et $\S
(\mathfrak{h})$ son alg\`{e}bre sym\'{e}trique.

Pour $x \in \mathfrak{h}$, on note $\sigma (x)$ l'unique
d\'{e}rivation de $\S (\mathfrak{h})$ qui prolonge l'endomorphisme $y
\to [x,y]$ de l'espace vectoriel $\mathfrak{h}$. Un id\'{e}al $J$ de
$\S (\mathfrak{h})$ est dit $\mathfrak{h}$-stable ou
$\mathfrak{h}$-invariant si $\sigma (x)(J) \subset J$ pour tout $x \in
\mathfrak{h}$. On d\'{e}signe par $\S (\mathfrak{h})^{\mathfrak{h}}$
l'ensemble des \'{e}l\'{e}ments $p$ de $\S (\mathfrak{h})$ tels que
$\sigma (x)(p) = 0$ pour tout $x \in \mathfrak{h}$.

\subsection Soient $\big( \U_{n}(\mathfrak{h})\big)_{n \geqslant 0}$
la filtration canonique de $\U (\mathfrak{h})$ et $\big(
\S^{n}(\mathfrak{h})\big)_{n \geqslant 0}$ la graduation canonique de
$\S (\mathfrak{h})$. On convient que $\U_{-1}(\mathfrak{h}) =
\S^{-1}(\mathfrak{h}) = \{ 0\}$. Si $a \in \U_{m}(\mathfrak{h})$ et $b
\in \U_{n}(\mathfrak{h})$, on a $ab - ba \in
\U_{m+n-1}(\mathfrak{h})$.

Soit $n \in \mathbb{N}$. Il existe un isomorphisme canonique d'espaces
vectoriels :
$$
j_{n} \colon \U_{n}(\mathfrak{h})/ \U_{n-1}(\mathfrak{h}) \to
\S^{n}(\mathfrak{h}) .
$$

On d\'{e}finit une structure de P-alg\`{e}bre sur $\S (\mathfrak{h})$
de la mani\`{e}re suivante : si $p \in \S^{m}(\mathfrak{h})$, $q \in
\S^{n}(\mathfrak{h})$, et si $\widetilde{p} \in \U_{m}(\mathfrak{h})$,
$\widetilde{q} \in \U_{n}(\mathfrak{h})$ v\'{e}rifient $p =
j_{m}(\widetilde{p})$, $q = j_{n}(\widetilde{q})$, on pose
$$
\{ p, q\} = j_{m+n-1}(\widetilde{p}\widetilde{q} -
\widetilde{q}\widetilde{p}). 
$$

La structure de P-alg\`{e}bre sur $\S (\mathfrak{h})$ ainsi
d\'{e}finie est dite canonique. Lorsque nous consid\`{e}rerons
$\S (\mathfrak{h})$ comme une P-alg\`{e}bre, ce sera toujours au moyen
de cette structure. 

\subsection Soient $x \in \mathfrak{h}$ et $p \in \S
(\mathfrak{h})$. On a :
$$
\{ x, p\} = \sigma (x)(p).
$$
On en d\'{e}duit qu'un id\'{e}al de $\S (\mathfrak{h})$ est un
P-id\'{e}al si et seulement s'il est $\mathfrak{h}$-stable. De
m\^{e}me, le centre de la P-alg\`{e}bre $\S (\mathfrak{h})$ est $\S
(\mathfrak{h})^{\mathfrak{h}}$. 

Soit $J$ un id\'{e}al $\mathfrak{h}$-stable de $\S
(\mathfrak{h})$. D'apr\`{e}s 2.4, $\S (\mathfrak{h}) /J$ est
canoniquement muni d'une structure de P-alg\`{e}bre.

\subsection Soient $n \in \mathbb{N}$ et $\mathfrak{h}$ l'alg\`{e}bre de
Lie de base $(x_{1}, y_{1}, \dots , x_{n}, y_{n},z)$ avec
$$
[x_{1}, y_{1}] = \cdots = [x_{n}, y_{n}] = z,
$$
les autres crochets \'{e}tant nuls ou s'en d\'{e}duisant par
antisym\'{e}trie. 

L'id\'{e}al $\S (\mathfrak{h})(z-1)$ est
$\mathfrak{h}$-stable. Les P-alg\`{e}bres $B_{n}$ et $\S
(\mathfrak{h})/ \S (\mathfrak{h})(z-1)$ sont isomorphes.

\section{Centre et semi-centre}

\chead[\sl Alg\`{e}bres de Poisson et alg\`{e}bres de Lie
r\'{e}solubles]{\sl 5~~Centre et semi-centre}

\subsection Dans le paragraphe 5, le corps $\Bbbk$ est suppos\'{e} 
alg\'{e}briquement clos. Les notations qui suivent seront
conserv\'{e}es dans toute la suite.

\subsection Soient $\mathfrak{h}$ une alg\`{e}bre de Lie et $V$ un
$\mathfrak{h}$-module.

Si $\lambda \in \mathfrak{h}^{*}$, on note $V_{\lambda}$ ou
$V_{\lambda}(\mathfrak{h})$ (resp. $V^{\lambda}$ ou
$V^{\lambda}(\mathfrak{h})$) le sous-espace de $V$ form\'{e} des
vecteurs $v$ qui v\'{e}rifient, pour tout $x \in \mathfrak{g}$, $x.v =
\lambda (x)v$ (resp. $\big( x - \lambda (x)\big)^{n}.v = 0$ d\`{e}s
que $n$ est assez grand). On a $V_{\lambda} \subset V^{\lambda}$, et
la somme des $V^{\lambda}$, pour $\lambda \in \mathfrak{h}^{*}$, est
directe. On dit que $\lambda$ est une forme lin\'{e}aire
distingu\'{e}e pour $V$ (resp. un poids de $V$) si $V^{\lambda} \ne \{
0\}$ (resp. $V^{\lambda} \ne \{ 0\}$).

On a le r\'{e}sultat suivant ([2], th\'{e}or\`{e}me 1.3.19).

\begin{th*}{\bf Proposition}On suppose $\mathfrak{h}$ nilpotente et
$V$ de dimension finie. 

{\em (i)} $V$ est la somme directe des $V^{\lambda}$, pour $\lambda
\in \mathfrak{g}^{*}$.

{\em (ii)} Pour tout $\lambda \in \mathfrak{g}^{*}$, l'espace
$V^{\lambda}$ est $\mathfrak{h}$-stable.

{\em (iii)} Si $x\in \mathfrak{h}$, notons $\rho (x)$ l'endomorphisme
de $V$ induit par $x$. Soit $\lambda \in \mathfrak{g}^{*}$. Il existe
une base $\mathcal{B}$ de $V^{\lambda}$ telle que, pour tout $x \in
\mathfrak{h}$, $\Mat \big( \rho (x) - \lambda (x) \id_{V},
\mathcal{B})$ soit triangulaire inf\'{e}rieure stricte. \end{th*}

\subsection Dans la suite de ce travail, $\mathfrak{g}$ est une
$\Bbbk$-alg\`{e}bre de Lie r\'{e}soluble de dimension finie et $Q$ un
id\'{e}al premier $\mathfrak{g}$-stable de $\S (\mathfrak{g})$.

On munit $B(Q,\mathfrak{g}) = B(Q) = \S (\mathfrak{g})/Q$ et $L(Q,
\mathfrak{g}) = L(Q) = \Fract B(Q)$ de leurs structures canoniques de
P-alg\`{e}bres (voir 2.4, 2.7 et 4.3). Si $x \in \mathfrak{g}$,
$\varepsilon_{Q}(x)$ est la P-d\'{e}rivation de $B(Q)$ ou de $L(Q)$
d\'{e}duite de $\ad x$. On note $Y(Q,\mathfrak{g}) = Y(Q)$
(resp. $D(Q,\mathfrak{g}) = D(Q)$) le centre de la P-alg\`{e}bre
$B(Q)$ (resp. $L(Q)$). Ainsi, $Y(Q)$ (resp. $D(Q)$) est l'ensemble des
$a \in B(Q)$ (resp. $a \in L(Q)$) qui v\'{e}rifient
$\varepsilon_{Q}(x)(a) = 0$ pour tout $x \in \mathfrak{g}$. On en
d\'{e}duit que $Y(Q)$ et $D(Q)$ sont des P-alg\`{e}bres commutatives. 

Soit $\lambda \in \mathfrak{g}^{*}$. Avec les notations de 5.2, un
\'{e}l\'{e}ment $a$ de $B(Q)_{\lambda}$ (resp. $L(Q)_{\lambda}$) est
appel\'{e} un semi-invariant de $B(Q)$ (resp. $L(Q)$). Si $a$ est non
nul, on dit que $\lambda$ est le poids de $a$ (cela n'entra\^{\i}ne
aucune confusion avec la terminologie de 5.2). Si $B(Q)_{\lambda} \ne
\{ 0\}$, on dit que $\lambda$ est distingu\'{e}e relativement \`{a}
$Q$. 

On note  $E(Q, \mathfrak{g}) = E(Q)$ l'ensemble des
semi-invariants non nuls de $B(Q)$. Il est clair que $E(Q)$ perme un
calcul de fractions dans $B(Q)$. 

On pose 
$$
SY (Q, \mathfrak{g}) = SY (Q) = \operatornamewithlimits{\textstyle
\sum}_{\lambda \in \mathfrak{g}^{*}} B(Q)_{\lambda} \ , \ SD(Q,
\mathfrak{g}) = SD( Q) = \operatornamewithlimits{\textstyle
  \sum}_{\lambda \in \mathfrak{g}^{*}} L(Q)_{\lambda},
$$
et on dit que $SY(Q)$ (resp. $SD(Q)$) est le semi-centre de $B(Q)$
(resp. $L(Q)$). 

Les r\'{e}sultats de 5.4 \`{a} 5.7 sont d\'{e}montr\'{e}s dans la
partie II de [9].

\begin{th}{\bf Lemme}{\em (i)} Soit $I$ un id\'{e}al non nul de
$B(Q)$. Il existe $\lambda \in \mathfrak{g}^{*}$ tel que $I \cap
B(Q)_{\lambda} \ne \{ 0\}$.

{\em (ii)} Soient $\lambda \in \mathfrak{g}^{*}$ et $a \in
L(Q)_{\lambda}$. Il existe $\mu \in \mathfrak{g}^{*}$, $b \in
B(Q)_{\lambda + \mu}$, et $c \in B(Q)_{\mu}$ tels que $a =
bc^{-1}$. \end{th} 

\subsection Posons
\begin{gather*}
\Lambda (Q,\mathfrak{g}) = \Lambda (Q) = \{ \lambda \in
\mathfrak{g}^{*}\, ; \, B(Q)_{\lambda} \ne \{ 0\}\}  , \\
 \Lambda'(Q, \mathfrak{g}) = \Lambda'(Q) = \{ \lambda \in
\mathfrak{g}^{*}\, ; \, L(Q)_{\lambda} \ne \{ 0\}\} .
\end{gather*}
D'apr\`{e}s 5.4, $\Lambda'(Q)$ est le sous-groupe additif de
$\mathfrak{g}^{*}$ engendr\'{e} par $\Lambda (Q)$. On a ainsi
$$
\operatornamewithlimits{\textstyle \bigcap}_{\lambda \in \Lambda (Q)}
\ker \lambda = \operatornamewithlimits{\textstyle \bigcap}_{\lambda
\in \Lambda'(Q)} \ker \lambda .
$$

\subsection Soit $\lambda \in \mathfrak{g}^{*}$ tel que $\lambda
([\mathfrak{g}, \mathfrak{g}]) = \{ 0\}$. Il existe un et un seul
automorphisme $\tau_{\lambda}$ de l'alg\`{e}bre $\S (\mathfrak{g})$ tel
que $\tau_{\lambda}(x) = x + \lambda (x)$ pour tout $x \in
\mathfrak{g}$. 

\begin{th*}{\bf Lemme}Soit $\lambda \in \mathfrak{g}^{*}$ v\'{e}rifiant
$\lambda ([\mathfrak{g}, \mathfrak{g}]) = \{ 0\}$, $\mathfrak{g}' =
\ker \lambda$, $Q$ un id\'{e}al premier $\mathfrak{g}$-stable de
$\S (\mathfrak{g})$, et $Q' = Q \cap \S (\mathfrak{g}')$. On a $Q =
\S (\mathfrak{g})Q'$ si et seulement si $\tau_{\lambda}(Q)= Q$. \end{th*}

\begin{th}{\bf Lemme}Soient $Q$ un id\'{e}al premier
$\mathfrak{g}$-stable de $\S (\mathfrak{g})$, $\lambda \in
\mathfrak{g}^{*} \backslash \{ 0\}$ v\'{e}rifiant $L(Q)_{\lambda} \ne
\{ 0\}$, $\mathfrak{g}'$ le noyau de $\lambda$, et $Q' =Q \cap
\S (\mathfrak{g}')$. On a $Q = \S (\mathfrak{g})Q'$. \end{th}

\begin{th}{\bf Lemme}Soient $\mathfrak{g}'$ un id\'{e}al de codimension
$1$ de $\mathfrak{g}$, $Q$ un id\'{e}al premier $\mathfrak{g}$-stable
de $\S (\mathfrak{g})$, et $Q' = Q \cap \S (\mathfrak{g}')$, de sorte
que $L(Q',\mathfrak{g}')$ s'identifie \`{a} un sous-corps de
$L(Q,\mathfrak{g})$. On suppose que $Q$ est l'id\'{e}al de
$\mathfrak{g}$ engendr\'{e} par $Q'$. Soient $\lambda \in
\mathfrak{g}^{*}$ et $a \in B(Q, \mathfrak{g})_{\lambda}$. Il existe
$b \in D(Q, \mathfrak{g})$ et $c \in B(Q , \mathfrak{g})_{\lambda}
\cap B(Q', \mathfrak{g}')$ tels que $a = bc$. \end{th}

\begin{Preuve}Soit $x \in \mathfrak{g} \backslash \mathfrak{g}'$. On
peut supposer que $a$ est non nul, et on d\'{e}signe par $u = x^{n}
u_{n} + x^{n-1}u_{n-1} + \cdots + u_{0}$ un repr\'{e}sentant de $a$
dans $\S (\mathfrak{g})$, avec $u_{0}, \dots , u_{n} \in
\S (\mathfrak{g}')$ et $u_{n} \notin Q'$.

Comme $[\mathfrak{g}, \mathfrak{g}] \subset \mathfrak{g}'$, si $y \in
\mathfrak{g}$ et $p \in \mathbb{N}^{*}$, il vient $\{ y, x^{p}\} \in
x^{p-1}S (\mathfrak{g}')$. On en d\'{e}duit que, modulo $Q$, on a :
\begin{align*}
\lambda (y) u & = \lambda (y) x^{n}u_{n} + \cdots + \lambda (y) u_{0}
\equiv \{ y, u\} \\
{} & \equiv x^{n} \{ y , u_{n}\} + x^{n-1} v_{n-1} + \cdots + v_{0},
\end{align*}
avec $v_{0}, \dots , v_{n-1} \in \S (\mathfrak{g}')$. Par
cons\'{e}quent, modulo $Q'$, on a $\{ y, u_{n}\} \equiv \lambda (y)
u_{n}$. Soit $c$ l'image de $u_{n}$ dans $B(Q',\mathfrak{g}')$. Comme
$u_{n} \notin Q'$, on a $c \ne 0$. Par suite, $ac^{-1}$ existe dans
$L(Q, \mathfrak{g})$, et il est imm\'{e}diat que $ac^{-1} \in D(Q,
\mathfrak{g})$, car $a$ et $c$ ont m\^{e}me poids
$\lambda$. \end{Preuve}

\begin{th}{\bf Lemme}Soient $\mathfrak{h}$ un id\'{e}al de
$\mathfrak{g}$, $Q$ un id\'{e}al premier $\mathfrak{g}$-stable de $\S
(\mathfrak{g})$, et $R = Q \cap \S (\mathfrak{h})$, de sorte que $B(R,
\mathfrak{h})$ s'identifie \`{a} une sous-alg\`{e}bre de $B(Q,
\mathfrak{g})$. Si $\mu$ est une forme lin\'{e}aire sur
$\mathfrak{h}$, distingu\'{e}e relativement \`{a} $R$, il existe
$\lambda \in \mathfrak{g}^{*}$, distingu\'{e}e relativement \`{a} $Q$,
et prolongeant $\mu$. \end{th}

\begin{Preuve}D'apr\`{e}s les hypoth\`{e}ses, $B(R, \mathfrak{h})$ est
un sous-$\mathfrak{g}$-module de $B(Q, \mathfrak{g})$. Comme
$\mathfrak{g}$ est r\'{e}soluble, $[\mathfrak{g}, \mathfrak{g}]$
op\`{e}re de mani\`{e}re localement nilpotente sur $B(Q,\mathfrak{g})$
et $B(R, \mathfrak{h})$.

Soient $a \in B(R, \mathfrak{h})_{\mu} \subset B(Q, \mathfrak{g})$, $x
\in \mathfrak{g}$ et $y \in \mathfrak{h}$. Il vient
\begin{align*}
\varepsilon_{Q}(y)\varepsilon_{Q}(x) (a) & = \varepsilon_{Q} ([y,x])
(a) + \varepsilon_{Q}(x) \varepsilon_{Q}(y) (a) \\
{} & = \mu ([y,x]) (a) + \mu (y) \varepsilon_{Q}(x)(a).
\end{align*}
D'apr\`{e}s ce qui pr\'{e}c\`{e}de, on a $\lambda ([y,x]) = 0$. Par
suite, $\varepsilon_{Q}(x)(a) \in B(R, \mathfrak{h})_{\mu}$. Comme
$\mathfrak{g}$ op\`{e}re de mani\`{e}re localement finie sur $B(R,
\mathfrak{h})_{\mu}$ et que $\mathfrak{g}$ est r\'{e}soluble, on a
obtenu le r\'{e}sultat. \end{Preuve}

\begin{th}{\bf Proposition}Soit $Q$ un id\'{e}al premier et
$\mathfrak{g}$-stable de $\S (\mathfrak{g})$. L'ensemble $SY(Q,
\mathfrak{g})$ est une {\em P}-sous-alg\`{e}bre commutative de $B(Q,
\mathfrak{g})$. \end{th}

\begin{Preuve}On raisonne par r\'{e}currence sur la dimension $n$ de
$\mathfrak{g}$, le r\'{e}sultat \'{e}tant clair pour $n \leqslant
1$. Si toute forme lin\'{e}aire sur $\mathfrak{g}$ distingu\'{e}e
relativement \`{a} $Q$ est nulle, on a $SY(Q, \mathfrak{g}) = Y(Q,
\mathfrak{g})$, et le r\'{e}sultat est \'{e}tabli. Supposons qu'il
existe $\lambda \in \mathfrak{g}^{*} \backslash \{ 0\}$ et $a \in B(Q,
\mathfrak{g})\backslash \{ 0\}$ v\'{e}rifiant $\varepsilon_{Q}(x)(a) =
\lambda (x)a$ pour tout $x \in \mathfrak{g}$. Soient $\mathfrak{g}' =
\ker \lambda$ et $Q' = Q \cap \S (\mathfrak{g}')$. On a $Q =\S
(\mathfrak{g}) Q'$ (5.7).

Soient $u_{1},u_{2} \in B(Q, \mathfrak{g})$ des
semi-invariants. D'apr\`{e}s 5.8, il existe des \'{e}l\'{e}ments
$v_{1}, v_{2} \in D(Q, \mathfrak{g})$ et $w_{1}, w_{2} \in SY(Q',
\mathfrak{g}')$ tels que $u_{1} = v_{1}w_{1}$ et $u_{2} =
u_{2}w_{2}$. Compte tenu de l'hypoth\`{e}se de r\'{e}currence, on a
$\{ w_{1}, w_{2}\} = 0$. Il vient alors $\{ u_{1}, u_{2}\} =
0$. \end{Preuve}

\begin{th}{\bf Lemme}Soient $\lambda$ une forme lin\'{e}aire non nulle
sur $\mathfrak{g}$, distingu\'{e}e relativement \`{a} $Q$,
$\mathfrak{g}'$ le noyau de $\lambda$, et $Q' = Q \cap \S
(\mathfrak{g}')$, de sorte que $B(Q', \mathfrak{g}')$ s'identifie
\`{a} une {\em P}-sous-alg\`{e}bre de $B(Q, \mathfrak{g})$.

{\em (i)} Soient $x \in \mathfrak{g} \backslash \mathfrak{g}'$,
$\delta$ la {\em P}-d\'{e}rivation de $B(Q', \mathfrak{g}')$ induite
par $x$, et $X$ une in\-d\'{e}\-ter\-mi\-n\'{e}e. Les {\em
P}-alg\`{e}bres $B(Q, \mathfrak{g})$ et $B(Q',
\mathfrak{g}')_{\delta}\{ X\}$, sont isomorphes.

{\em (ii)} On a $SY (Q, \mathfrak{g}) \subset B(Q',
\mathfrak{g}')$. \end{th}

\begin{Preuve}(i) Soit $\widetilde{x}$ l'image de $x$ dans
$B(Q, \mathfrak{g})$. Avec des notations \'{e}videntes notons $\theta$
l'application
$$ 
B(Q', \mathfrak{g}')_{\delta}\{ X\} \to B(Q, \mathfrak{g}) \ , \
\operatornamewithlimits{\textstyle \sum} X^{n} a_{n} \to
\operatornamewithlimits{\textstyle \sum} \widetilde{x}^{n} a_{n}.
$$
Il est imm\'{e}diat que $\theta$ est un homomorphisme surjectif de
P-alg\`{e}bres. Comme $Q = \S (\mathfrak{g}) Q'$ (5.6), on voit que
$\theta$ est bijectif.

(ii) Notons encore $\delta$ la d\'{e}rivation de $L(Q',
\mathfrak{g}')$ induite par $\ad x$. Si $S$ est l'ensemble des
\'{e}l\'{e}ments non nuls de $B(Q', \mathfrak{g}')$, on a $B(Q,
\mathfrak{g})_{S} = L(Q', \mathfrak{g}')_{\delta} \{ X\}$. Prouvons
que $\delta$ n'est pas une d\'{e}rivation int\'{e}rieure de $L(Q',
\mathfrak{g}')$. En effet, soit $e \in B(Q, \mathfrak{g})_{\lambda}
\backslash \{ 0\}$. Si $\delta = d_{u}$, avec $u \in L(Q',
\mathfrak{g}')$, on obtient $\delta (e) = \{ u, e\} = 0$ (par
d\'{e}finition de $\mathfrak{g}'$). Or, $\delta (e) = \{ X, e\} =
\lambda (x) e \ne 0$. D'o\`{u} l'assertion.

Soit $a \in B(Q, \mathfrak{g})\backslash \{ 0\}$ un semi-invariant de
poids $\mu$. D'apr\`{e}s 5.8, il s'\'{e}crit $a = bc$, avec $b \in
D(Q, \mathfrak{g})$ et $c \in B(Q, \mathfrak{g})_{\mu} \cap B(Q',
\mathfrak{g}')$. D'apr\`{e}s ce qui pr\'{e}c\`{e}de et 2.9,\,(ii), on
a $b \in D(Q', \mathfrak{g}')$. Ainsi, $a \in L(Q', \mathfrak{g}')
\cap B(Q, \mathfrak{g}) = B(Q', \mathfrak{g}')$. \end{Preuve}

\subsection Dans la suite, si $Q$ est un id\'{e}al premier
$\mathfrak{g}$-stable de $\S (\mathfrak{g})$, on pose :
$$
\widehat{\mathfrak{g}}_{Q} = \operatornamewithlimits{\textstyle
\bigcap}_{\lambda \in \Lambda (Q)} \ker \lambda \ , \ \widehat{Q} = Q
\cap S(\widehat{\mathfrak{g}}_{Q}) .
$$

Comme $[\mathfrak{g}, \mathfrak{g}] \subset
\widehat{\mathfrak{g}}_{Q}$, si $(x_{1}, \dots , x_{r})$ est une base
d'un suppl\'{e}mentaire de $\widehat{\mathfrak{g}}_{Q}$ dans
$\mathfrak{g}$, alors $\widehat{\mathfrak{g}}_{Q} + \Bbbk x_{1} +
\cdots + \Bbbk x_{i}$ est un id\'{e}al de $\mathfrak{g}$ pour $1
\leqslant i \leqslant r$.

\begin{th}{\bf Th\'{e}or\`{e}me}Soient $\mathfrak{g}$ une alg\`{e}bre
de Lie r\'{e}soluble et $Q$ un id\'{e}al premier
$\mathfrak{g}$-stable de $\S (\mathfrak{g})$.

{\em (i)} On a $SY(Q, \mathfrak{g}) \subset SY (\widehat{Q},
\widehat{\mathfrak{g}}_{Q}) = Y (\widehat{Q},
\widehat{\mathfrak{g}}_{Q})$.

{\em (ii)} Soient $(x_{1}, \dots , x_{r})$ une base d'un
suppl\'{e}mentaire de $\widehat{\mathfrak{g}}_{Q}$ dans
$\mathfrak{g}$. En notant $\delta_{i}$, $1 \leqslant i \leqslant r$,
la d\'{e}rivation de $B(Q, \mathfrak{g})$ induite par $x_{i}$, alors
les {\em P}-alg\`{e}bres $B(Q, \mathfrak{g})$ et $\big( \cdots \big(
B(\widehat{Q}, \widehat{\mathfrak{g}}_{Q})_{\delta_{1}}\{ X_{1}\}\big)
\cdots \big)_{\delta_{r}} \{ X_{r}\}$ sont isomorphes. \end{th}

\begin{Preuve}(i) Le fait que $SY (Q, \mathfrak{g}) \subset SY
(\widehat{Q}, \widehat{\mathfrak{g}}_{Q})$ se d\'{e}duit facilement de
5.11,(ii).

Supposons qu'il existe une forme lin\'{e}aire $\mu$ non nulle sur
$\widehat{\mathfrak{g}}_{Q}$ et distingu\'{e}e relativement \`{a}
$\widehat{Q}$. D'apr\`{e}s 5.9, il existe $\lambda \in
\mathfrak{g}^{*}$ distingu\'{e}e relativement \`{a} $Q$, et
prolongeant $\mu$. On a alors $\widehat{\mathfrak{g}}_{Q} \subset \ker
\lambda$, donc $\widehat{\mathfrak{g}}_{Q} \subset
\widehat{\mathfrak{g}}_{Q} \cap \ker \lambda = \ker \mu$. Ainsi, $\mu
= 0$. Contradiction. On a donc bien $SY (\widehat{Q},
\widehat{\mathfrak{g}}_{Q}) = Y(\widehat{Q},
\widehat{\mathfrak{g}}_{Q})$.

(ii) On raisonne par r\'{e}currence sur la dimension de $\mathfrak{g}$,
le cas o\`{u} $\dim \mathfrak{g} = 0$ \'{e}tant clair. On peut
supposer que $\dim (\mathfrak{g} / \widehat{\mathfrak{g}}_{Q}) = r >
0$. 

Soit $(\lambda_{1}, \dots , \lambda_{r})$ une base de l'orthogonal de
$\widehat{\mathfrak{g}}_{Q}$ dans $\mathfrak{g}^{*}$ telle que l'on
ait $\lambda_{i}(x_{j}) = \delta_{ij}$ pour $1 \leqslant i,j \leqslant
r$. Soient $\mathfrak{h} = \ker \lambda_{r} =
\widehat{\mathfrak{g}}_{Q} + \Bbbk x_{1} + \cdots + \Bbbk x_{r-1}$ et $Q' = Q
\cap \S (\mathfrak{h})$. D'apr\`{e}s 5.11,\,(i), les alg\`{e}bres
$B(Q, \mathfrak{g})$ et $B(Q', \mathfrak{h})_{\delta_{r}}\{ X_{r}\}$
sont isomorphes. On termine alors facilement d'apr\`{e}s
l'hypoth\`{e}se de r\'{e}currence car, compte tenu de 5.9, on a
$\widehat{\mathfrak{h}}_{Q'} =
\widehat{\mathfrak{g}}_{Q}$. \end{Preuve}

\section{Formes lin\'{e}aires distingu\'{e}es}

\chead[\sl Alg\`{e}bres de Poisson et alg\`{e}bres de Lie
r\'{e}solubles]{\sl 6~~Formes lin\'{e}aires distingu\'{e}es}

\subsection On conserve les hypoth\`{e}ses et notations du paragraphe
5.

Soit $V$ un $\mathfrak{g}$-module de dimension finie $n$. Il
poss\`{e}de une suite de Jordan-H\"{o}lder, c'est-\`{a}-dire qu'il
existe une suite de $\mathfrak{g}$-modules
$$
\{ 0\} = V_{0} \subset V_{1} \subset \cdots \subset V_{n} = V
$$ 
v\'{e}rifiant $\dim V_{i} = i$ pour $0 \leqslant i \leqslant
n$. D'autre part, si $(V'_{i})_{0 \leqslant i \leqslant n}$ est une
autre suite de Jordan-H\"{o}lder de $V$, il existe une permutation
$\sigma$ de l'ensemble $\{ 1, 2, \dots , n\}$ telle que les
$\mathfrak{g}$-modules $V'_{i}/ V'_{i-1}$ et $V_{\sigma (i)}/
V_{\sigma (i) -1}$ soient isomorphes pour $1 \leqslant i \leqslant
n$.

Chaque $\mathfrak{g}$-module $V_{i}/ V_{i-1}$ \'{e}tant de dimension
$1$ s'identifie \`{a} une forme lin\'{e}aire $\mu_{i}$ sur
$\mathfrak{g}$. On note $\mathscr{J}(V,\mathfrak{g})$ ou
$\mathscr{J}(V)$ l'ensemble $\{ \mu_{1}, \dots , \mu_{n}\}$. 

Soient $W$ un sous-module de $V$ et $(W_{i})_{0 \leqslant i \leqslant
r}$ une suite de Jordan-H\"{o}lder de $W$. Il existe des sous-modules 
$$
W = V_{r} \subset \cdots \subset V_{n} = V
$$ 
tels que $(V_{i}/W)_{r \leqslant i \leqslant n}$ soit une suite de
Jordan-H\"{ol}der de $V/W$. Alors
$$
\{ 0\} = W_{0} \subset \cdots \subset W_{r} = W = V_{r} \subset \cdots
\subset V_{n} = V
$$
est une suite de Jordan-H\"{o}lder de $V$. On en d\'{e}duit que l'on a
:
\begin{equation}\label{3}
\mathscr{J}(W) \subset \mathscr{J}(V) \ \text{ et } \ \mathscr{J}(V/W)
\subset \mathscr{J}(V).
\end{equation}

\subsection Consid\'{e}rons $\mathfrak{g}$ comme un
$\mathfrak{g}$-module au moyen de la repr\'{e}sentation adjointe, et
fixons une suite 
$$
\{ 0\} = \mathfrak{g}_{0} \subset \mathfrak{g}_{1} \subset \cdots
\subset \mathfrak{g}_{m} = \mathfrak{g}
$$
d'id\'{e}aux de $\mathfrak{g}$ formant une suite de Jordan-H\"{o}lder
de $\mathfrak{g}$ (on a donc $\dim \mathfrak{g} = m$). Soit
$\lambda_{i} \in \mathfrak{g}^{*}$ la forme lin\'{e}aire correspondant
\`{a} $\mathfrak{g}_{i}/ \mathfrak{g}_{i-1}$ pour $1 \leqslant i
\leqslant m$. 

Soit $\mathcal{B} = (y_{1}, \dots , y_{m})$ une base de $\mathfrak{g}$
telle que $(y_{1}, \dots , y_{i})$ soit une base de $\mathfrak{g}_{i}$
pour $1 \leqslant i \leqslant n$. On a ainsi
$$
[x,y_{i}] \in \lambda_{i}(x) y_{i} + \mathfrak{g}_{i-1}
$$
pour $1 \leqslant i \leqslant m$ et tout $x \in \mathfrak{g}$.

Si $\nu = (\nu_{1}, \dots , \nu_{m}) \in \mathbb{N}^{m}$, on pose
$|\nu | = \nu_{1} + \cdots + \nu_{m}$, et on note
$y^{\nu}$ l'\'{e}l\'{e}ment de $\S (\mathfrak{g})$ d\'{e}fini par
$$
y^{\nu} = y_{1}^{\nu_{1}} y_{2}^{\nu_{2}} \cdots y_{m}^{\nu_{m}}.
$$

On d\'{e}finit un ordre total sur l'ensemble des mon\^{o}mes $y^{\nu}$
en convenant, avec des notations \'{e}videntes que $y^{\nu'} <
y^{\nu}$ si l'une des conditions suivantes est r\'{e}alis\'{e}e :

$\bullet$ ou $|\nu' | < |\nu |$

$\bullet$ ou $|\nu' | = |\nu |$ et il existe $i \in \{ 1, 2, \dots ,
m\}$ tel que
$$
\nu'_{m} = \nu_{m}, \dots , \nu'_{i+1} = \nu_{i+1} \ \text{ et } \
\nu'_{i} < \nu_{i}.
$$

Si $x \in \mathfrak{g}$ et $\nu = (\nu_{1}, \dots , \nu_{m}) \in
\mathbb{N}^{m}$, on a facilement 
$$
\{ x, y^{\nu}\} = \big( \nu_{1} \lambda_{1}(x) + \cdots + \nu_{m}
\lambda_{m} (x)\big) y^{\nu} + v,
$$
o\`{u} $v$ est une combinaison lin\'{e}aire de termes $y^{\nu'}$ avec
$|\nu'| = |\nu|$ et $y^{\nu'} < y^{\nu}$.

Si $n \in \mathbb{N}$, avec les notations de 4.2, il est alors
imm\'{e}diat que
$$
\mathscr{J} \big( \S^{n}(\mathfrak{g})\big) \subset \mathbb{N}
\lambda_{1} + \cdots + \mathbb{N} \lambda_{n}.
$$

Compte tenu de \eqref{3}, on obtient alors facilement :

\begin{th*}{\bf Lemme}Soient $Q$ un id\'{e}al $\mathfrak{g}$-stable de
$\S (\mathfrak{g})$ et $V$ un $\mathfrak{g}$-module de dimension
finie de $\S (\mathfrak{g})/ Q$. On a :
$$
\mathscr{J}(V) \subset \mathbb{N} \lambda_{1} + \cdots + \mathbb{N}
\lambda_{m} .
$$
En particulier :
$$ 
\Lambda (Q) \subset \mathbb{N} \lambda_{1} + \cdots + \mathbb{N}
\lambda_{m}.
$$ \end{th*}

\subsection Conservons les notations $\lambda_{1}, \dots ,
\lambda_{m}$ pr\'{e}c\'{e}dentes, et supposons que le plus grand
id\'{e}al nilpotent $\mathfrak{n}$ de $\mathfrak{g}$ soit commutatif. 

Soient $\lambda \in \mathfrak{g}^{*} \backslash \{ 0\}$, $x \in
\mathfrak{g} \backslash \ker \lambda$, et $y \in \mathfrak{g}
\backslash \mathfrak{n}$. On a
$$
\big( \ad x - \lambda (x) \id_{\mathfrak{g}}\big) (y) = [x,y] -
\lambda (x) y.
$$ 
Comme $\lambda (x) \ne 0$ et $[x,y] \in [\mathfrak{g},
\mathfrak{g}] \subset \mathfrak{n}$, on d\'{e}duit de ceci que $y
\notin \mathfrak{g}^{\lambda}$ (notation de 5.2). Il en r\'{e}sulte
que, pour tout poids $\nu$ de $\mathfrak{g}$, on a $\mathfrak{g}^{\nu}
\subset \mathfrak{n}$.

Comme $\mathfrak{n}$ est commutatif, on peut consid\'{e}rer
$\mathfrak{n}$ comme un $\ad (\mathfrak{g}
/\mathfrak{n})$-module. L'alg\`{e}bre de Lie $\mathfrak{g} /
\mathfrak{n}$ \'{e}tant nilpotente, on peut appliquer 5.2. Il existe
donc des formes lin\'{e}aires deux \`{a} deux distinctes $\mu_{1},
\dots , \mu_{s}$ sur $\mathfrak{g}$ telles que
$$
\mathfrak{n} = \mathfrak{n}^{\mu_{1}} \oplus \cdots \oplus
\mathfrak{n}^{\mu_{s}} .
$$ 

D'autre part, pour $1 \leqslant i \leqslant s$, on a $[\mathfrak{g},
\mathfrak{n}^{\mu_{i}}] \subset \mathfrak{n}^{\mu_{i}}$. Enfin, il
existe une base $\mathcal{B}_{i}$ de $\mathfrak{n}^{\mu_{i}}$ telle
que la matrice de $\big( \ad x - \mu_{i}(x) \id_{\mathfrak{g}}\big) |
\mathfrak{n}^{\mu_{i}}$ soit triangulaire inf\'{e}rieure stricte.

Soit $\mathcal{B}' = (e_{q+1}, \dots , e_{m})$ la base de
$\mathfrak{n}$ obtenue par r\'{e}union des bases
$\mathcal{B}_{i}$. Compl\'{e}tons pour obtenir une base $\mathcal{B}
= (e_{1}, \dots , e_{m})$ de $\mathfrak{g}$. Comme $[\mathfrak{g},
\mathfrak{g}] \subset \mathfrak{n}$, pour tout $x \in \mathfrak{g}$,
la matrice $\Mat (\ad x, \mathcal{B})$ est de la forme
$$
\begin{pmatrix}
0 & 0 \\
A(x) & T(x) \end{pmatrix},
$$
o\`{u} $A(x)$ a $\dim \mathfrak{n}$ lignes et $\dim
(\mathfrak{g}/\mathfrak{n})$ colonnes, et o\`{u} $T(x)$ est
carr\'{e}e, triangulaire inf\'{e}rieure, et a $\dim \mathfrak{n}$
lignes. 

On d\'{e}duit en particulier de ceci que tout \'{e}l\'{e}ment non nul
de $\mathscr{J}(\mathfrak{g})$ est l'un des $\mu_{i}$ et que l'on a
le r\'{e}sultat suivant.

\begin{th*}{\bf Lemme}On suppose que le plus grand id\'{e}al nilpotent
de $\mathfrak{g}$ est commutatif. Si $\lambda$ est un \'{e}l\'{e}ment
non nul de $\mathscr{J}(\mathfrak{g})$, il existe $y \in \mathfrak{g}
\backslash \{ 0\}$ tel que
$$
[x,y] = \lambda (x)y
$$ 
pour tout $x \in \mathfrak{g}$. \end{th*}

\section{Une classe d'alg\`{e}bres de Poisson}

\chead[\sl Alg\`{e}bres de Poisson et alg\`{e}bres de Lie
r\'{e}solubles]{\sl 7~~Une classe d'alg\`{e}bres de Poisson}

\subsection Dans ce paragraphe, on ne suppose plus que $\Bbbk$ est
alg\'{e}briquement clos. Les notations qui suivent seront
utilis\'{e}es dans toute la suite.

Soient $V$ un $\Bbbk$-espace vectoriel de dimension finie $n$ et $\S
(V)$ l'alg\`{e}bre sy\-m\'{e}\-tri\-que de $V$. On d\'{e}signe par
$\omega$ une forme bilin\'{e}aire altern\'{e}e sur $V$ et par $G$ un
sous-groupe libre de type fini du dual $V^{*}$ de $V$. Soit $p = \rg
(G)$ le rang de $G$. Si $g \in G$, $\lambda_{g}$ est la forme
lin\'{e}aire sur $V$ correspondant \`{a} $g$. On note $\Bbbk [G]$
l'alg\`{e}bre du groupe $G$, $V^{\omega}$ le noyau de $\omega$,
$V^{G}$ l'orthogonal de $G$ dans $V$, c'est-\`{a}-dire 
$$
V^{G} = \operatornamewithlimits{\textstyle \bigcap}_{g \in G} \ker
\lambda_{g}, 
$$
et $V^{G \omega}$ le noyau de la restriction de $\omega$ \`{a} $V^{G}
\times V^{G}$. 
 
Soient $X_{1}, \dots , X_{n}, Y_{1}, \dots , Y_{p}$ des
ind\'{e}termin\'{e}es. Les alg\`{e}bres associatives $\S (V)
\otimes_{\Bbbk} \Bbbk [G]$ et $\Bbbk [X_{1}, \dots , X_{n}, Y_{1},
Y_{1}^{-1}, \dots , Y_{p}, Y_{p}^{-1}]$ sont isomorphes, et les
unit\'{e}s de $\S (V) \otimes_{\Bbbk} \Bbbk [G]$ sont les
\'{e}l\'{e}ments de la forme $\mu g$, avec $\mu \in \Bbbk \backslash
\{ 0\}$ et $g \in G$.

\subsection Il existe une et une seule structure de P-alg\`{e}bre sur
$\S (V) \otimes_{\Bbbk} \Bbbk [G]$ telle que, pour $v,w \in V$ et $g,h
\in G$, on ait :
$$
\{ v, w\} = \omega (v,w) \ , \ \{ g,h \} = 0 \ , \ \{ g,v\} =
\lambda_{g}(v)g. 
$$

La P-alg\`{e}bre ainsi obtenue est not\'{e}e $\mathscr{B}_{\Bbbk}(V,
\omega ,G)$ ou $\mathscr{B}(V, \omega ,G)$. On note $\S_{\omega}(V)$
la P-sous-alg\`{e}bre de $\mathscr{B}(V,\omega ,G)$ engendr\'{e}e par
$V$. Si l'on consid\`{e}re $\S (V^{\omega})$ comme une P-alg\`{e}bre
commutative, il est clair que les P-alg\`{e}bres $\S_{\omega}(V)$ et
$B_{r} \otimes_{\Bbbk} S(V^{\omega})$, o\`{u} $2r$ est le rang de
$\omega$, sont isomorphes.

\subsection Le r\'{e}sultat suivant est de d\'{e}monstration
imm\'{e}diate.

\begin{th*}{\bf Proposition}Soient $B$ une {\em P}-alg\`{e}bre et $U$
le groupe de ses unit\'{e}s. On d\'{e}signe par $\chi \colon V \to
B$ une application $\Bbbk$-lin\'{e}aire, et par $\psi \colon G \to U$ un
homomorphisme de groupes v\'{e}rifiant les conditions suivantes :

{\em (i)} $\{ \chi (v), \chi (v')\} = \omega (v,v')$ pour tous $v,v'
\in V$.

{\em (ii)} $\{ \psi (g), \chi (v)\} = \lambda_{g}(v) \psi (g)$ pour
tout $g \in G$ et tout $v \in V$.

{\em (iii)} $\{ \psi (g), \psi (g')\} = 0$ pour tous $g,g' \in G$.

Alors il existe un unique homomorphisme de {\em P}-alg\`{e}bres
$\theta \colon \mathscr{B}(V, \omega , G) \to B$ tel que $\theta | V =
\chi$ et $\theta | G = \psi$. \end{th*}

\subsection Soit $\lambda \in V^{*}$. Il existe une et une seule
d\'{e}rivation $D_{\lambda}$ de l'alg\`{e}bre associative $\S (V)$
telle que $D_{\lambda}(v) = \lambda (v)$ pour tout $v \in V$. Cette
d\'{e}rivation est localement nilpotente. On peut donc d\'{e}finir
$\exp D_{\lambda}$. Si $g \in G$, on note $\varphi_{g}$ pour $\exp
D_{\lambda_{g}}$. 

\begin{th*}{Lemme}{\em (i)} Soient $g \in G$ et $r \in \S_{\omega}(V)$.
Dire que $\varphi_{g}(r) = r$ signifie que $\{ g, r\} = 0$.

{\em (ii)} On suppose que $\{ \lambda_{g}|V^{\omega}\, ; \, g \in G\}$
engendre l'espace vectoriel $(V^{\omega})^{*}$. Soit $r \in \S
(V^{\omega}) \backslash \Bbbk$. Il existe $g \in G$ tel que
$\varphi_{g}(r) \ne r$. \end{th*}

\begin{Preuve}(i) Soit $V' = \ker \lambda_{g}$. Si $r \in \S (V')$, il
est imm\'{e}diat que $\varphi_{g}(r) = r$ et $\{ g, r\} =
0$. Supposons $r \in \S_{\omega} (V) \backslash \S (V')$. Il existe $v
\in V$ tel que $\lambda_{g}(v) = 1$, et on peut \'{e}crire
$$
r = v^{n} r_{n} + v^{n-1}r_{n-1} + \cdots + r_{0},
$$
avec $n \geqslant 1$, $r_{0}, \dots , r_{n} \in \S (V')$ et $r_{n} \ne
0$. Il vient :
\begin{gather*}
\{ g,r\} = g(n v^{n-1}r_{n-1} + (n-1) v^{n-2} r_{n-2} + \cdots +
r_{1}), \\ 
\varphi_{g}(r) = (v+1)^{n}r_{n} + (v+1)^{n-1} r_{n-1} + \cdots +
(v+)r_{1} + r_{0}.
\end{gather*}
On a donc $\varphi_{g}(r) \ne r$ et $\{ g, r\} \ne 0$.

(ii) Soient $g_{1}, \dots , g_{n} \in G$ tels que
$(\lambda_{g_{1}}|V^{\omega}, \dots \lambda_{g_{n}}|V_{\omega})$
soit une base de $(V^{\omega})^{*}$, et soit $(v_{1}, \dots ,
v_{n})$ sa base duale dans $V^{\omega}$. 

Pour $\nu_{1}, \dots , \nu_{n} \in \mathbb{N}$ et $1 \leqslant i
\leqslant n$, on a :
$$
\varphi_{g_{i}}(v_{1}^{\nu_{1}} \cdots v_{n}^{\nu_{n}}) =
v_{1}^{\nu_{1}} \cdots v_{i-1}^{\nu_{i-1}} (v_{i} +1)^{\nu_{i}}
\nu_{i+1}^{\nu_{i+1}} \cdots v_{n}^{\nu_{n}}.
$$
L'assertion en d\'{e}coule facilement. \end{Preuve}

\begin{th}{\bf Proposition}Les conditions suivantes sont
\'{e}quivalentes :

{\em (i)} La {\em P}-alg\`{e}bre $\mathscr{B}(V, \omega , G)$ est
simple.

{\em (ii)} L'ensemble des $r \in \S (V^{\omega })$ tels que
$\varphi_{g}(r) = r$ pour tout $g \in G$ est r\'{e}duit \`{a} $\Bbbk$.

{\em (iii)} L'ensemble $\{ \lambda_{g}| V^{\omega}\, ; g \in G\}$
engendre l'espace vectoriel $(V^{\omega})^{*}$, ce qui signifie que
$V^{G} \cap V^{\omega} = \{ 0\}$.

{\em (iv)} Les seuls id\'{e}aux $J$ de $\S (V^{\omega})$ v\'{e}rifiant
$\varphi_{g}(J) \subset J$ pour tout $g \in G$ sont $\{ 0\}$ et $\S
(V^{\omega})$. \end{th}

\begin{Preuve}(i) $\Rightarrow$ (ii) Soit $r \in \S (V^{\omega})
\backslash \Bbbk$ v\'{e}rifiant $\varphi_{g}(r) = r$ pour tout $g \in
G$. D'apr\`{e}s 7.4,\,(i), $r$ est un \'{e}l\'{e}ment central non
inversible de $\mathscr{B}$. Par suite, $\mathscr{B}$ n'est pas une
P-alg\`{e}bre simple.

(ii) $\Rightarrow$ (iii) Si $r \in V^{\omega}$, alors $\varphi_{g}(r)
= r + \lambda_{g}(r)$ pour tout $g \in G$. L'implication est donc
claire. 

(iii) $\Rightarrow$ (ii) R\'{e}sulte de 7.4,\,(ii).

(iii) $\Rightarrow$ (iv) Soit $J$ unid\'{e}al non nul de $\S
(V^{\omega})$ tel que $\varphi_{g}(J) \subset J$ pour tout $g \in
G$. Supposons $J \ne \S (V^{\omega})$, et soit $r \in J$ de degr\'{e}
minimal $m$ parmi les \'{e}l\'{e}ments non nuls de $J$. On a $m
\geqslant 1$. D'apr\`{e}s 7.4,\,(ii), il existe $g \in G$ tel que
$\varphi_{g}(r) \ne r$. Le degr\'{e} de $r - \varphi_{g}(r)$ \'{e}tant
strictement inf\'{e}rieur \`{a} celui d e$r$, c'est absurde.

(iv) $\Rightarrow$ (ii) Soit $r \in \S (V^{\omega}) \backslash \{ 0\}$
v\'{e}rifiant $\varphi_{g}(r) = r$ pour tout $g \in G$. Alors $J = \S
(V^{\omega}) r$ est un id\'{e}al de $\S (V^{\omega})$ tel que
$\varphi_{g}(J) \subset J$ pour tout $G \in G$. On en d\'{e}duit que
$J = \S (V^{\omega})$, donc $r$ est inversible dans $\S
(V^{\omega})$. D'o\`{u} $r \in \Bbbk \backslash \{ 0\}$.

(iii) $\Rightarrow$ (i) Soient $v \in V$, $u \in \S_{\omega}(V)
\backslash \{ 0\}$, et $g \in G$. Il vient :
$$
\{ v, ug\} = \big( \{v,u\}- \lambda_{g}(v)u \big) g.
$$

Si $\lambda_{g}(v) \ne 0$, on a donc $\{ v, ug\} \ne 0$ car, ou $\{ v,
u\} = 0$, ou le degr\'{e} de $\{ v, u\}$  est strictement
inf\'{e}rieur \`{a} celui de $u$. 

Soit $J$ un P-id\'{e}al non nul de $\mathscr{B}$. Tout \'{e}l\'{e}ment
non nul $r$ de $J$ s'\'{e}crit
$$
r = r_{1} g_{1} + \cdots + r_{s}g_{s},
$$
avec $r_{1}, \dots , r_{s} \in \S_{\omega}(V)$ et $g_{1}, \dots , g_{s}
\in G$. Choisissons $r$ de mani\`{e}re que $s$ soit minimal et que le
degr\'{e} de $r_{1}$ soit minimal. On peut supposer que $g_{1}$ est
l'\'{e}l\'{e}ment neutre de $G$. 

On a $r_{1} \in \S (V^{\omega})$. En effet, sinon, il existe $v \in V$
tel que $\{ v,r_{1}\} \ne 0$ et tel que le degr\'{e} de $\{ v, r_{1}\}$
soit strictement inf\'{e}rieur \`{a} celui de $r_{1}$. Cela contredit
le choix de $r$.

Supposons $r_{2} \ne 0$. Alors $g_{2}$ est distinct de
l'\'{e}l\'{e}ment neutre de $G$ (d'apr\`{e}s la minimalit\'{e} de
$s$). Soit $v \in V \backslash \ker \lambda_{g_{2}}$. Il vient :
$$
\{ v , r\} = \operatornamewithlimits{\textstyle \sum}_{i=2}^{s} \big(
\{ v, r_{i}\} - \lambda_{g_{i}}(v) r_{i}\big) g_{i}.
$$
Comme $\lambda_{g_{2}}(v) \ne 0$, ceci contredit la minimalit\'{e} de
$s$ d'apr\`{e}s ce qui pr\'{e}c\`{e}de.

On a donc $r_{2} = 0$ et, on a ainsi prouv\'{e} que $J \cap
\S (V^{\omega}) \ne \{ 0\}$. 

Reprenons les notations de la preuve de 7.4,\,(ii). Pour $\nu_{1},
\dots , \nu_{n} \in \mathbb{N}$ et $1 \leqslant i \leqslant n$, on
obtient :
$$
g_{i}^{-1}\{ g_{i}, v_{1}^{\nu_{1}} \cdots v_{n}^{\nu_{n}}\} = \nu_{i}
v_{1}^{\nu_{1}} \cdots v_{i-1}^{\nu_{i-1}} v_{i}^{\nu_{i}-1}
v_{i+1}^{\nu_{i+1}} \cdots v_{n}^{\nu_{n}}.
$$
On en d\'{e}duit que les seuls id\'{e}aux de $\S (V^{\omega})$ stables
par les applications $r \to g^{-1}\{g,r\}$, $g \in G$, sont $\{ 0\}$
et $\S (V^{\omega})$. Par cons\'{e}quent, $J \cap \S (V^{\omega}) = \S
(V^{\omega})$ et $J = \mathscr{B}$. \end{Preuve}

\begin{th}{\bf Proposition}Soit $\mathscr{B} = \mathscr{B}(V, \omega , G)$
une {\em P}-alg\`{e}bre simple.

{\em (i)} Le centralisateur $C$ de $\Bbbk [G]$ dans $\mathscr{B}$ est
$S_{\omega}(V^{G}) \otimes_{\Bbbk} \Bbbk [G]$.

{\em (ii)} Le centre $D$ de $C$ est $S_{\omega}(V^{G\omega}) \otimes_{\Bbbk}
\Bbbk [G]$. \end{th}

\begin{Preuve}(i) Il est imm\'{e}diat que $S_{\omega}(V^{G})
\otimes_{\Bbbk} \Bbbk [G] \subset C$.

Soient $W$ un suppl\'{e}mentaire de $V^{G}$ dans $V$. Il existe des
\'{e}l\'{e}ments $g_{1}, \dots , g_{n}$ de $G$ tels que $(g_{1}, \dots
, g_{n})$ soit une base du dual $W^{*}$ de $W$. Soit $(v_{1}, \dots ,
v_{n})$ la base de $W$ duale de la pr\'{e}c\'{e}dente. Un
\'{e}l\'{e}ment $a$ de $\mathscr{B}$ s'\'{e}crit
$$
a = \operatornamewithlimits{\textstyle \sum}_{i_{1}, \dots , i_{n}}
v_{1}^{i_{1}} \cdots v_{n}^{i_{n}} a_{i_{1}, \dots ,i_{n}}, 
$$
avec $i_{1}, \dots , i_{n} \in \mathbb{N}$ et $a_{i_{1}, \dots ,
i_{n}} \in S_{\omega} (V^{G}) \otimes_{\Bbbk} \Bbbk [G]$.

On a, par exemple :
$$
\{ g_{1}, a\} = \operatornamewithlimits{\textstyle \sum}_{i_{1}, \dots
, i_{n}} i_{1} v_{1}^{i_{1}-1} v_{2}^{i_{2}} \cdots v_{n}^{i_{n}}
g_{1} a_{i_{1}, \dots  ,i_{n}}.
$$
On voit donc que $a \in C$ si et seulement si $a = a_{0, \dots , 0}$,
soit $a \in S_{\omega} (V^{G}) \otimes_{\Bbbk} \Bbbk [G]$.

(ii) Il existe une base $(u_{1}, v_{1}, \dots , u_{r}, v_{r}\}$ d'un
suppl\'{e}mentaire de $V^{G\omega}$ dans $V^{G}$ telle que
$$
\{ u_{1}, u_{j}\} = \{ v_{i}, v_{j}\} = 0 \ , \ \{ u_{i}, v_{j}\} =
\omega_{ij} 
$$
pour $1 \leqslant i,j \leqslant r$. En utilisant ceci, on obtient
facilement le r\'{e}sultat. \end{Preuve}

\subsection Si $A$ est une $\Bbbk$-alg\`{e}bre associative, on note
$\mathbf{d} (A)$ sa dimension de Gelfand-Kirillov (voir
[3]). Rappelons que, si $A$ est commutative, int\`{e}gre, et de type
fini, alors $\mathbf{d}(A)$ est le degr\'{e} de transcendance sur
$\Bbbk$ du corps des fractions de $A$.

\begin{th}{\bf Proposition}Soit $\mathscr{B} = \mathscr{B}(V, \omega , G)$
une {\em P}-alg\`{e}bre simple comme en {\em 7.2}. On note $\rg
(G)$ le rang de $G$, et on conserve les notations $C$ et $D$ de {\em
7.6}. Alors :

{\em (i)} $\mathbf{d} \big( \mathscr{B}(V,\omega , G)\big) = \dim
V + \rg (G)$.

{\em (ii)} $\mathbf{d} \big( \Bbbk [G] \big) = \rg (G)$.

{\em (iii)} $\mathbf{d} (C) = \dim V^{G} + \rg (G)$.

{\em (iv)} $\mathbf{d}(D) = \dim V^{G\omega} + \rg (G)$.

En particulier, les entiers $\rg (G)$, $\dim V$, $\dim V^{G}$ et
$\dim V^{G\omega}$, ne d\'{e}pendent que de $\mathscr{B}$ et non de sa
pr\'{e}sentation sous la forme $\mathscr{B}(V, \omega ,G)$. \end{th} 

\begin{Preuve}Si $(v_{1}, \dots , v_{n})$ est une base de $V$ et
$(g_{1}, \dots , g_{r})$ une base du $\mathbb{Z}$-module $G$, on a
$\mathscr{B} = \Bbbk [v_{1}, \dots , v_{n}, g_{1}, g_{1}^{-1}, \dots ,
g_{r}, g_{r}^{-1}]$ en tant qu'alg\`{e}bre associative. On en
d\'{e}duit imm\'{e}diatement (i). Les autres assertions sont analogues
compte tenu de 7.6. \end{Preuve}

\subsection Soient $x_{1}, \dots , x_{n}, y_{1}, \dots , y_{n}$ des
g\'{e}n\'{e}rateurs de la P-alg\`{e}bre $B_{n}$ (voir 3.1)
v\'{e}rifiant
$$
\{ x_{i}, x_{j}\} = \{ y_{i}, y_{j}\} = 0 \ , \ \{ x_{i}, y_{j}\} =
\delta_{ij} 
$$
pour $1 \leqslant i,j \leqslant j$.

On note $B'_{n}$ la $P$-alg\`{e}bre localis\'{e}e de $B_{n}$ par les
puissances de $x = x_{1} \cdots x_{n}$. 

\begin{th}{\bf Proposition}Soit $\mathscr{B} = \mathscr{B}(V, \omega , G)$
une $P$-alg\`{e}bre simple comme en {\em 7.2}.

{\em (i)} Si $m$ est le rang de $G$ et $2 \ell$ celui de $\omega$,
$\mathscr{B}$ est une sous-{\em P}-alg\`{e}bre de la {\em P}-alg\`{e}bre
$B_{\ell} \otimes_{k} B'_{m}$.

{\em (ii)} Il existe une $\Bbbk$-alg\`{e}bre de Lie r\'{e}soluble
$\mathfrak{g}$ et un 
id\'{e}al premier $\mathfrak{g}$-stable $Q$ de $\S (\mathfrak{g})$ tels
que $\mathscr{B}$ soit isomorphe \`{a} la {\em P}-alg\`{e}bre $\big(
\S (\mathfrak{g}) /Q\big)_{E}$, o\`{u} $E$ est l'ensemble des
semi-invariants non nuls de $\S (\mathfrak{g})/ Q$. \end{th}

\begin{Preuve}(i) Fixons une base  $(x_{1}, y_{1}, \dots , x_{\ell},
y_{\ell}, s_{1}, \dots , s_{t})$ de $V$ telle que $(s_{1}, \dots ,
s_{t})$ soit une base du noyau de $\omega$ et telle que
$$
\{ x_{i}, x_{j}\} = \{ y_{i}, y_{j}\} = 0 \ , \ \{ x_{i}, y_{j}\} =
\delta_{ij}
$$
pour $1 \leqslant i,j \leqslant n$. 

Soit $\{ g_{1}, \dots , g_{m}\}$ une base du $\mathbb{Z}$-module
$G$. Il existe des scalaires $\lambda_{ji}, \mu_{ji}, \nu_{jp}$ tels
que
$$ \{ g_{j},x_{i}\} = \lambda_{ji} g_{j} \ , \ \{ g_{j}, y_{i}\} =
\mu_{ji}g_{j} \ , \ \{ g_{j}, s_{p}\} = \nu_{ji} g_{j}
$$
pour $1 \leqslant i \leqslant \ell$, $1 \leqslant j \leqslant m$, $1
\leqslant p \leqslant t$.

Ecrivons 
$$
B_{\ell} = \Bbbk [X_{1}, Y_{1}, \dots , X_{\ell}, Y_{\ell}] \ , \ B'_{m} =
\Bbbk [ Z_{1}, Z_{1}^{-1}, T_{1}, \dots , Z_{m}, Z_{m}^{-1}, T_{m}], 
$$
avec 
\begin{gather*}
\{ X_{i}, X_{j}\} = \{ Y_{i}, Y_{j}\} = 0 \ , \ \{ X_{i}, Y_{j}\} =
\delta_{ij} ,\\
\{ Z_{p}, Z_{q}\} = \{ T_{p}, T_{q}\} = 0 \ , \ \{ Z_{p}, T_{q}\} =
\delta_{pq}  
\end{gather*}
pour $1 \leqslant i,j \leqslant \ell$ et $1 \leqslant p,q \leqslant
m$. 

D\'{e}finissons une application lin\'{e}aire 
$$
\chi \colon V \to \operatornamewithlimits{\textstyle
\sum}_{i=1}^{\ell} \Bbbk X_{i} + \operatornamewithlimits{\textstyle
\sum}_{i=1}^{\ell} \Bbbk X_{i} + \operatornamewithlimits{\textstyle
\sum}_{i=1}^{m} \Bbbk T_{i}
$$
en posant 
$$
\chi (x_{i}) = X_{i} + \operatornamewithlimits{\textstyle
  \sum}_{j=1}^{m} \lambda_{ji} T_{j} \ , \ \chi (y_{i}) = Y_{i} +
\operatornamewithlimits{\textstyle \sum}_{j=1}^{m} \mu_{ji} T_{j} \ ,
\ \chi (s_{p}) = \operatornamewithlimits{\textstyle \sum}_{j=1}^{m}
\nu_{jp} T_{j}
$$
pour $1 \leqslant i \leqslant \ell$ et $1 \leqslant p \leqslant t$.

De m\^{e}me, soit $\psi \colon G \to \Bbbk [Z_{1}, Z_{1}^{-1}, \dots ,
Z_{m}, Z_{m}^{-1}]$ l'application $g_{i} \to Z_{i}$ pour $1 \leqslant
i \leqslant m$.

Compte tenu de 7.3, ces applications se prolongent en un
homomorphisme de P-alg\`{e}bres $\mathscr{B} \to B_{\ell} \otimes_{\Bbbk}
B'_{m}$. On obtient alors l'assertion car, $\mathscr{B}$ \'{e}tant
simple, cet homomorphisme est injectif. 

(ii) Conservons le notations pr\'{e}c\'{e}dentes. Soit $\mathfrak{g}$
l'alg\`{e}bre de Lie de dimension $2 \ell + t + m+1$ d\'{e}finie par
$$ 
\mathfrak{g} = \Bbbk w \oplus V \oplus \operatornamewithlimits{\textstyle
\sum}_{i=1}^{m} \Bbbk g_{i},
$$
avec 
$$
[x_{i}, x_{j}] = [y_{i}, y_{j}] = 0 \ , \ [x_{i}, y_{j}] =
\delta_{ij}w \ ,\  [g_{p}, g_{q}] = 0 \ , \ [g_{p}, v] = 
\lambda_{g_{p}}(v) g_{p}
$$
pour $1 \leqslant i,j \leqslant \ell$, $1 \leqslant p,q \leqslant m$,
et $v \in V$.

Soit $Q = S (\mathfrak{g})(w -1)$. Alors $Q$ est un id\'{e}al premier
$\mathfrak{g}$-stable de $S (\mathfrak{g})$. Si l'on note
$\widetilde{g}_{i}$ l'image de $g_{i}$ dans $\S (\mathfrak{g})/Q$, alors
$E$ est le semi-groupe engendr\'{e} par les $g_{i}$, $1 \leqslant i
\leqslant m$. On voit alors que $\mathscr{B}$ et $\big(
\S (\mathfrak{g})/ Q\big)_{E}$ sont isomorphes d'apr\`{e}s 7.3 et la
simplicit\'{e} de $\mathscr{B}$. \end{Preuve}

\section{Un cas particulier}

\chead[\sl Alg\`{e}bres de Poisson et alg\`{e}bres de Lie
r\'{e}solubles]{\sl 8~~Un cas particulier}

\subsection Dans toute la suite, le corps $\Bbbk$ est suppos\'{e}
alg\'{e}briquement clos. On reprend la notation $B_{n}$ de 3.1. On
dira que des \'{e}l\'{e}ments $x_{1}, y_{1}, \dots , x_{n}, y_{n}$ de
$B_{n}$ sont des g\'{e}n\'{e}rateurs de $B_{n}$ s'ils engendrent
l'alg\`{e}bre associative sous-jacente \`{a} $B_{n}$ et si
$$
\{ x_{i}, x_{j}\} = \{ y_{i}, y_{j}\} = 0 \ , \ \{ x_{i}, y_{j}\} = \delta_{ij}
$$
pour $1 \leqslant i,j \leqslant n$.

\begin{th}{\bf Th\'{e}or\`{e}me}Soient $\mathfrak{h}$ une alg\`{e}bre
de Lie r\'{e}soluble, $\mathfrak{g}$ un id\'{e}al de $\mathfrak{h}$,
et $\mathfrak{s}$ une sous-alg\`{e}bre de Lie de $\mathfrak{h}$
op\'{e}rant sur $\mathfrak{g}$ de mani\`{e}re semi-simple. On
d\'{e}signe par $Q$ un id\'{e}al premier de $S( \mathfrak{g})$ tel que
$\{ x, Q\} \subset Q$ pour tout $x \in \mathfrak{h}$, et on fait
op\'{e}rer $\mathfrak{h}$ sur $B(Q, \mathfrak{g})$ au moyen de la
repr\'{e}sentation d\'{e}duite de la repr\'{e}sentation adjointe. On
suppose que $SY (Q , \mathfrak{g}) = Y(Q, \mathfrak{g})$. Il existe $e
\in Y(Q, \mathfrak{g}) \backslash \{ 0\}$ et $n \in \mathbb{N}$
v\'{e}rifiant les conditions suivantes :

{\em (i)} L'alg\`{e}bre $Y(Q, \mathfrak{g})_{e}$ est de type fini.

{\em (ii)} $e$ est un vecteur propre pour l'action de $\mathfrak{h}$
dans $B(Q, \mathfrak{g})$.

{\em (iii)} Les {\em P}-alg\`{e}bres $B(Q,\mathfrak{g})_{e}$ et $Y(Q,
\mathfrak{g})_{e} \otimes_{\Bbbk} B_{n}$ sont isomorphes. 

{\em (iv)} On peut choisir des g\'{e}n\'{e}rateurs $x_{1}, y_{1},
\dots , x_{n}, y_{n}$ de $B_{n} \subset B(Q,\mathfrak{g})_{e}$ qui
soient des vecteurs propres pour l'action de $\mathfrak{s}$ dans $B(Q,
\mathfrak{g})_{e}$. \end{th}

\begin{Preuve}On raisonne par r\'{e}currence sur la dimension de
$\mathfrak{g}$, le r\'{e}sultat \'{e}tant \'{e}vident pour $\dim
\mathfrak{g} \leqslant 1$.

Soit $\mathfrak{g}'$ un id\'{e}al de $\mathfrak{h}$, contenu dans
$\mathfrak{g}$, et v\'{e}rifiant $\dim (\mathfrak{g}/ \mathfrak{g}') =
1$. On pose $Q' = Q \cap \S (\mathfrak{g}')$. L'image de $\S
(\mathfrak{g}')$ dans $B(Q, \mathfrak{g})$ s'identifie \`{a} $B(Q',
\mathfrak{g}')$.

D'apr\`{e}s 5.9, on a $SY (Q', \mathfrak{g}') = Y(Q',
\mathfrak{g}')$. Compte tenu de l'hypoth\`{e}se de
r\'{e}\-cur\-ren\-ce, il existe $n \in \mathbb{N}$ et $e' \in Y(Q',
\mathfrak{g}')$ vecteur propre pour l'action de $\mathfrak{h}$,
tels que
$$ 
B(Q', \mathfrak{g}')_{e'} = Y(Q', \mathfrak{g}')_{e'} \otimes_{\Bbbk}
B_{n},
$$
et v\'{e}rifiant les conditions suivantes :

(i) L'alg\`{e}bre $Y(Q', \mathfrak{g}')_{e'}$ est de type fini.

(ii) La P-sous-alg\`{e}bre $B_{n}$ de $B(Q', \mathfrak{g}')_{e'}$ a
des g\'{e}n\'{e}rateurs $x_{1}, y_{1}, \dots , x_{n}, y_{n}$ qui sont
vecteurs propres pour l'action de $\mathfrak{s}$ dans $B(Q',
\mathfrak{g}')_{e'}$.

Soient $z \in \mathfrak{g} \backslash \mathfrak{g}'$, $\widetilde{z}$
son image dans $B(Q, \mathfrak{g})$, et $\delta$ la P-d\'{e}rivation
de $B(Q', \mathfrak{g}')_{e'}$ induite par $\ad z$. D'apr\`{e}s les
hypoth\`{e}ses, on peut supposer qu'il existe une forme lin\'{e}aire
$\theta$ sur $\mathfrak{s}$ telle que $[t,z] = \theta (t)z$ pour tout
$t \in \mathfrak{s}$. La P-sous-alg\`{e}bre $Y(Q',
\mathfrak{g}')_{e'}$ de $B(Q', \mathfrak{g}')_{e'}$ est stable par
$\delta$, et la restriction de $\delta$ \`{a} $Y(Q',
\mathfrak{g}')_{e'}$ est localement nilpotente puisque $SY (Q,
\mathfrak{g}) = Y(Q, \mathfrak{g})$.

a) Supposons $\delta \big( Y(Q', \mathfrak{g}')\big) = \{
0\}$, donc $e' \in Y(Q, \mathfrak{g})$. D'apr\`{e}s 3.4,\,(iv), il
existe $b \in B(Q', \mathfrak{g}')_{e'}$ tel que $\delta =
d_{b}$. Posons $x = \widetilde{z} - b$. La P-alg\`{e}bre $B(Q,
\mathfrak{g})_{e'}$ est engendr\'{e}e par $x$ et $B(Q',
\mathfrak{g}')_{e'}$, et $x$ commute \`{a} $B(Q',
\mathfrak{g}')_{e'}$. On a donc :
$$ 
B(Q, \mathfrak{g})_{e'} = \big( B(Q',\mathfrak{g}')_{e'}\big) \{
x\} = \big( Y(Q', \mathfrak{g}')_{e'}\big) \{ x\} \otimes_{\Bbbk} B_{n}.
$$ 
Il est imm\'{e}diat que $Y(Q, \mathfrak{g})_{e'} = \big( Y(Q',
\mathfrak{g}')_{e'}\big) \{ x\}$, et que $Y(Q, \mathfrak{g})_{e'}$ est
une alg\`{e}bre de type fini. On a obtenu le r\'{e}sultat.

b) Dans la suite, on suppose $\delta \big( Y(Q', \mathfrak{g}')\big)$
non r\'{e}duit \`{a} $\{ 0\}$. 
 
Soient $u \in Y(Q', \mathfrak{g}')$, $h \in \mathfrak{h}$, et
$\Delta$ la P-d\'{e}rivation de $B(Q', \mathfrak{g}')$ induite par
$\ad h$. Il vient :
$$
\Delta \big( \delta (u)\big) = \Delta ( \{ \widetilde{z}, u\} ) = \{
\Delta (\widetilde{z}) , u\} + \{ \widetilde{z}, \Delta (u)\} .
$$ 
Comme $[\mathfrak{h}, \mathfrak{g}'] \subset \mathfrak{g}'$,
$[\mathfrak{h}, \mathfrak{g}] \subset \mathfrak{g}$, que $\mathfrak{g}
= \mathfrak{g}' + \Bbbk z$, et que $Y(Q', \mathfrak{g}')$ est le
centre de $B(Q', \mathfrak{g}')$, on obtient $\Delta \big( \delta
(u)\big) \in \delta \big( Y(Q', \mathfrak{g}')\big)$. Ainsi, $\delta
\big( Y(Q', \mathfrak{g}')\big)$ est un $\mathfrak{h}$-sous-module de
$Y(Q', \mathfrak{g}')$.

L'alg\`{e}bre de Lie $\mathfrak{h}$ \'{e}tant r\'{e}soluble, il
existe $v \in \delta \big( Y(Q', \mathfrak{g}')\big) \backslash \{
0\}$ vecteur propre pour l'action de $\mathfrak{h}$. Comme $v$ est
aussi vecteur propre pour l'action de $\mathfrak{g}$, on a $v \in SY
(Q, \mathfrak{g}) = Y(Q, \mathfrak{g})$. Ainsi, $\delta (v) = 0$.

Notons $V$ l'ensemble des $u \in Y(Q', \mathfrak{g}')$ qui
v\'{e}rifient $\delta (u) \in \Bbbk v$. Si $u \in V$ et $h \in
\mathfrak{h}$, avec les notations pr\'{e}c\'{e}dentes, il vient :
$$
\Delta \big( \delta (u)\big) = \Delta (\{ \widetilde{z}, u\} ) = \{
\Delta (\widetilde{z}) , u\} + \{ \widetilde{z}, \Delta (u)\} = 
\{ \Delta(\widetilde{z}), u\} + \delta \big( \Delta (u)\big). 
$$
Pour les m\^{e}mes raisons que dans l'alin\'{e}a pr\'{e}c\'{e}dent, on
voit alors que $V$ est un $\mathfrak{h}$-sous-module de $Y(Q',
\mathfrak{g}')$. D'autre part, $\mathfrak{s}$ op\`{e}re de mani\`{e}re
localement finie dans $B(Q, \mathfrak{g})$. On en d\'{e}duit qu'il
existe $u \in V$, vecteur propre pour l'action de $\mathfrak{s}$, et
tel que $\delta (u) = v$.

c) Soient $v \in \delta \big( Y(Q', \mathfrak{g}')\big) \backslash
\{0\}$ comme dans le point b) et $u \in Y(Q', \mathfrak{g}')$ vecteur
propre pour l'action de $\mathfrak{s}$ et v\'{e}rifiant $\delta (u) =
v$. Posant $y = uv^{-1}$, on a $\delta (y) = 1$ car $\delta (v) = 0$,
et $y$ est vecteur propre pour l'action de $\mathfrak{s}$. Notons
enfin $e = e'v$. Il vient $e \in Y(Q', \mathfrak{g}')$.

D'apr\`{e}s 3.4,\,(iv), il existe $b \in B(Q', \mathfrak{g}')_{e'}$
tel que $\delta' = \delta - d_{b}$ v\'{e}rifie :
$$ 
\delta '(1 \otimes B_{n}) = 0 \ , \ \delta' | Y(Q',
\mathfrak{g}')_{e'} \otimes 1 = \delta | Y(Q', \mathfrak{g}')_{e'}
\otimes 1.
$$ 
Posons $x = \widetilde{z} - b$. On a $\{ x, y\} = 1$ et, les seuls
P-id\'{e}aux de $B_{1}$ \'{e}tant $\{ 0\}$ et $B_{1}$ (3.4), la
P-sous-alg\`{e}bre de $B(Q, \mathfrak{g})_{e}$ engendr\'{e}e par $x$
et $y$ est isomorphe \`{a} $B_{1}$. 

Prouvons que l'on peut choisir $b$ de sorte que $x$ soit vecteur
propre pour l'action de $\mathfrak{s}$.

Si $\lambda \in \mathfrak{s}^{*}$, notons $W_{\lambda}$ l'ensemble
des $u \in B(Q', \mathfrak{g'})_{e'}$ qui v\'{e}rifient
$\varepsilon_{Q'}(t)(u) = \lambda (t) u$ pour tout $t \in
\mathfrak{s}$. D'apr\`{e}s les hypoth\`{e}ses, il vient 
$$ 
B(Q',\mathfrak{g}')_{e'} = \operatornamewithlimits{\textstyle
\bigoplus}_{\lambda \in \mathfrak{s}^{*}} W_{\lambda}.
$$

Reprenons les notations $x_{1}, y_{1}, \dots , x_{n}, y_{n}$ du
d\'{e}but de la preuve. Il existe $\mu_{i}, \nu_{i} \in
\mathfrak{s}^{*}$ tels que $x_{i} \in W_{\lambda_{i}}$ et $y_{i} \in
W_{\nu_{i}}$ pour $1 \leqslant i \leqslant n$. Rappelons d'autre part
que $\widetilde{z} \in W_{\theta}$. 

Ecrivons $b = b_{1} + \cdots  + b_{r}$, o\`{u} $b_{i} \in W_{\varphi_{i}}
\backslash \{ 0\}$, et o\`{u} $\varphi_{1}, \dots , \varphi_{r}$ sont
des formes lin\'{e}aires deux \`{a} deux distinctes sur
$\mathfrak{s}$. Pour $1 \leqslant i \leqslant n$, on a 
$$
0 = \{ x, x_{i}\} = \{ \widetilde{z}, x_{i}\} -
\operatornamewithlimits{\textstyle \sum}_{j=1}^{r}\{ b_{j}, x_{i}\} \
, \ 0 = \{ x, y_{i}\} = \{ \widetilde{z}, y_{i}\} -
\operatornamewithlimits{\textstyle \sum}_{j=1}^{r} \{ b_{j}, y_{i}\} .
$$
D'autre part :
$$
\{ \widetilde{z}, x_{i}\} \in W_{\theta + \lambda_{i}} \ , \ \{
\widetilde{z}, y_{i}\} \in W_{\theta + \nu_{i}} \ , \ \{ b_{j},
x_{i}\} \in W_{\varphi_{j} + \lambda_{i}} \ , \ \{ b_{j}, y_{i}\} \in
W_{\varphi_{j} + \nu_{i}}.
$$
Il est donc imm\'{e}diat que l'on peut supposer que
$\varepsilon_{Q'}(t)(b) = \theta (t) b$ pour tout $t \in
\mathfrak{s}$ ; on s'est ramen\'{e} au cas o\`{u} $x$ est
vecteur propre pour l'action de $\mathfrak{s}$. 

Ainsi, la P-sous-alg\`{e}bre de $B(Q,\mathfrak{g})_{e}$ engendr\'{e}e
par $x$ et $y$ est isomorphe \`{a} $B_{1}$, et $x,y$ sont vecteurs
propres pour l'action de $\mathfrak{s}$.

d) D'apr\`{e}s 3.9, $\big( Y(Q', \mathfrak{g}')_{e}\big)_{\delta}\{ X\}$
et $\big( Y(Q', \mathfrak{g}')_{e} / Y(Q', \mathfrak{g}')_{e}y \big)
\otimes_{\Bbbk} B_{1}$ sont des P-alg\`{e}bres isomorphes.

Soient $X,Y$ des ind\'{e}termin\'{e}es. Consid\'{e}rons l'application
$$
\tau \colon \big( B(Q', \mathfrak{g}')_{e}\big)_{\delta'}\{ X\} \to
\big( B(Q', \mathfrak{g}')_{e}\big)_{\delta}\{ Y\}
$$ 
d\'{e}finie, pour $p \in \mathbb{N}$, $u \in B(Q',
\mathfrak{g}')_{e}$, et $c \in B_{n}$, par
$$
\tau [(u \otimes c)X^{p}] = (u \otimes c) (Y-b)^{p}.
$$
Il est clair que $\tau$ est un isomorphisme d'alg\`{e}bres
associatives. D'autre part :
$$
\tau (\{ X, u \otimes c\}) = \tau (\delta' (u)\otimes c) + \tau
(u \otimes \delta'(c)\} = \delta (u) \otimes c.
$$

De m\^{e}me, dans $\big( B(Q', \mathfrak{g}')_{e}\big)_{\delta}\{
Y\}$, on a :
$$
\{ Y-b, u \otimes c\} = \{ Y-b,u\} \otimes c + u \otimes \{ Y-b, c\} =
\delta (u) \otimes c.
$$
Par suite, $\tau$ est un isomorphisme de P-alg\`{e}bres. 

Il r\'{e}sulte de tout ceci que les P-alg\`{e}bres
\begin{gather*} 
\big( B(Q', \mathfrak{g}')_{e}\big)_{\delta}\{ X\} \ , \ \big(
B(Q', \mathfrak{g}')_{e}\big)_{\delta'}\{ X\}  , \\ 
\big( Y(Q', \mathfrak{g}')_{e}\big)_{\delta'} \{ X\} \otimes_{\Bbbk}
B_{n} \ , \ \big( Y(Q', \mathfrak{g}')_{e}\big)_{\delta} \{ X\}
\otimes_{\Bbbk} B_{n}
\end{gather*}
sont isomorphes. 

On en d\'{e}duit que les P-alg\`{e}bres 
$$
\big( B(Q', \mathfrak{g}')_{e}\big)_{\delta}\{ X\} \text{ et } \ \big(
Y(Q', \mathfrak{g}')_{e} / 
Y(Q', \mathfrak{g}')_{e}y) \big) \otimes_{\Bbbk} B_{n+1} 
$$
sont isomorphes. 

Il est clair que $B(Q, \mathfrak{g})_{e}$ est un quotient non nul de
$\big( B(Q', \mathfrak{g'})_{e}\big)_{\delta} \{ X\}$. Compte tenu de
3.4, $B(Q, \mathfrak{g})_{e}$ est de la forme $Z \otimes B_{n+1}$,
o\`{u} $Z$ est un quotient non nul de $Y(Q', \mathfrak{g}')_{e} /
Y(Q', \mathfrak{g}')_{e}y$. On en d\'{e}duit que $Z$ est contenu dans
le centre de $B(Q, \mathfrak{g})_{e}$ et est de type fini. Les
relations de 3.1 montrent alors que $Z$ est le centre de $B(Q,
\mathfrak{g})_{e}$. \end{Preuve}

\begin{th}{\bf Corollaire}Soient $\mathfrak{g}$ une alg\`{e}bre de Lie
nilpotente, $Q$ un id\'{e}al premier $\mathfrak{g}$-stable de $\S
(\mathfrak{g})$, $E$ l'ensemble des \'{e}l\'{e}ments centraux non nuls
de la {\em P}-alg\`{e}bre $B(Q) = \S (\mathfrak{g})/Q$, et $D(Q)$ le
centre de $\Fract B(Q)$. Il existe un entier $n \in \mathbb{N}$ tel
les {\em P}-alg\`{e}bres $B(Q)_{E}$ et $D(Q) \otimes_{\Bbbk} B_{n} =
B_{n} \big( D(Q)\big)$ soient isomorphes.  \end{th}

\begin{Preuve}Avec les notations de 8.2, on a
$D(Q) = Y(Q)_{E} = \big( Y(Q)_{e}\big)_{E}$. D'o\`{u} le
r\'{e}sultat. \end{Preuve}

\subsection Soient $\mathfrak{g}$ une alg\`{e}bre de Lie et $\Gamma$
son groupe alg\'{e}brique adjoint. Le groupe $\Gamma$ op\`{e}re
naturellement dans $\mathfrak{g}$. On le fait op\'{e}rer dans
$\mathfrak{g}^{*}$ en convenant que
$$
(\gamma .f)(x) = f (\gamma^{-1} .x)
$$
pour $x \in \mathfrak{g}$, $f \in \mathfrak{g}^{*}$, et $\gamma \in
\Gamma$. 

\begin{th*}{\bf Corollaire}Soient $\mathfrak{g}$ une alg\`{e}bre de Lie
nilpotente, $\Gamma$ son groupe adjoint, et $Q$ un id\'{e}al premier
$\mathfrak{g}$-stable de $\S (\mathfrak{g})$. Les conditions suivantes
sont \'{e}quivalentes :

{\em (i)} $Y (Q) = \Bbbk$.

{\em (ii)} Il existe $n \in \mathbb{N}$ tel que $B(Q) = B_{n}$.

{\em (iii)} $Q$ est un id\'{e}al $\mathfrak{g}$-stable maximal.

{\em (iv)} $Q$ est l'id\'{e}al de $\S (\mathfrak{g})$ associ\'{e}
\`{a} une $\Gamma$-orbite dans $\mathfrak{g}^{*}$. \end{th*}

\begin{Preuve}L'alg\`{e}bre de Lie $\mathfrak{g}$ \'{e}tant
nilpotente, elle op\`{e}re de mani\`{e}re localement nilpotente dans
$B(Q)$. On a donc $SY(Q) = Y(Q)$. D'o\`{u} l'implication (i)
$\Rightarrow$ (ii) d'apr\`{e}s 8.2. On obtient aussi (ii)
$\Rightarrow$ (i) d'apr\`{e}s le lemme 3.4.

D'apr\`{e}s le hypoth\`{e}ses, on a $\gamma (u) \in Q$ pour tout $u
\in Q$. Le groupe $\Gamma$ op\`{e}re naturellement dans $B(Q)$, et le
centre de $B(Q)$ est l'ensemble des points fixes de $\Gamma$ dans
$B(Q)$. 

Compte tenu de ces remarques et de [2], propositions 4.8.5 et
7.3.2,\,(ii), on obtient les \'{e}quivalences (i) $\Leftrightarrow$
(iii) $\Leftrightarrow$ (iv). \end{Preuve}

\section{Un second cas particulier}

\chead[\sl Alg\`{e}bres de Poisson et alg\`{e}bres de Lie
r\'{e}solubles]{\sl 9~~Un second cas particulier}

\subsection Dans ce paragraphe, $\mathfrak{g}$ est une alg\`{e}bre de
Lie r\'{e}soluble dont le plus grand id\'{e}al nilpotent
$\mathfrak{n}$ est commutatif. On d\'{e}signe par $Q$ un id\'{e}al
premier $\mathfrak{g}$-stable de $\S (\mathfrak{g})$ tel que $Q \cap
\mathfrak{g} = \{ 0\}$. On identifiera souvent un \'{e}l\'{e}ment de
$\mathfrak{g}$ \`{a} son image dans $B(Q) = \S (\mathfrak{g})/ Q$.

On note $\lambda_{1}, \dots , \lambda_{s}$ les \'{e}l\'{e}ments deux
\`{a} deux distincts de $\mathscr{J}(\mathfrak{g}) \backslash \{ 0\}$
(voir 6.1 et 6.2). On a donc $\mathfrak{n} \subset \ker \lambda_{1}
\cap \cdots \ker \lambda_{n}$. R\'{e}ciproquement, si $x \in \ker
\lambda_{i}$ pour $1 \leqslant i \leqslant s$, il r\'{e}sulte de  [10],
19.5.7 que $x \in \mathfrak{n}$. 

D'apr\`{e}s 6.3, si $1 \leqslant i \leqslant s$, il existe $e_{i} \in
\mathfrak{g}$ tel que $[x, e_{i}] = \lambda_{i}(x) e_{i}$ pour tout $x
\in \mathfrak{g}$. Comme $Q \cap \mathfrak{g} = \{ 0\}$ et que $Q$ est
premier, on d\'{e}duit de 6.1 que
$$
\Lambda (Q) = \mathbb{N} \lambda_{1} + \cdots + \mathbb{N} \lambda_{s}.
$$

Posons $e = e_{1} \cdots e_{s}$. L'id\'{e}al $Q$ \'{e}tant premier et
v\'{e}rifiant $Q \cap \mathfrak{g} = \{ 0\}$, on a $e^{n} \notin Q$
pour tout $n \in \mathbb{N}$. On peut donc consid\'{e}rer le
$\mathfrak{g}$-module (ou la P-alg\`{e}bre) $B(Q)_{e}$ et, l'action de
$\mathfrak{g}$ dans $B(Q)_{e}$ \'{e}tant localement finie, il
r\'{e}sulte de 6.1 que l'ensemble des poids des semi-invariants de
$B(Q)_{e}$ est
$$
G = \mathbb{Z} \lambda_{1} + \cdots + \mathbb{Z} \lambda_{s}.
$$

Le corps $\Bbbk$ \'{e}tant de caract\'{e}ristique nulle, $G$ est sans
torsion, donc est un groupe libre de rang fini $\ell$. Fixons
$\mu_{1}, \dots , \mu_{\ell} \in \mathfrak{g}^{*}$ tels que
$$
G = \mathbb{Z} \mu_{1} + \cdots + \mathbb{Z} \mu_{\ell},
$$
et soient $\alpha_{1}, \dots , \alpha_{\ell} \in B(Q)_{e}$ des
semi-invariants non nuls, de poids respectifs $\mu_{1}, \dots ,
\mu_{\ell}$, et qui sont des mon\^{o}mes en $e_{1}, e_{1}^{-1}, \dots
, e_{s}, e_{s}^{-1}$.

Si $n_{1}, \dots , n_{\ell} \in \mathbb{Z}$, $\alpha_{1}^{n_{1}}
\cdots \alpha_{\ell}^{n_{\ell}}$ est semi-invariant de poids $n_{1}
\mu_{1} + \cdots + n_{\ell} \mu_{\ell}$. Il en r\'{e}sulte que, si
$\mathscr{E} = \{ \alpha_{1}, \alpha_{1}^{-1}, \dots , \alpha_{\ell},
\alpha_{\ell}^{-1}\} $, des mon\^{o}mes distincts en des
\'{e}l\'{e}ments de $\mathscr{E}$ ont des poids distincts. Par suite,
le groupe multiplicatif $\Gamma$ engendr\'{e} par les \'{e}l\'{e}ments
de $\mathscr{E}$ est isomorphe \`{a} $G$. De m\^{e}me, la famille des
mon\^{o}mes distincts en les \'{e}l\'{e}ments de $\mathscr{E}$ est
libre sur $\Bbbk$ (voir 5.2). On en d\'{e}duit que la sous-alg\`{e}bre
de $B(Q)_{e}$ engendr\'{e}e par $e_{1}, \dots , e_{n}, e^{-1}$
s'identifie \`{a} l'alg\`{e}bre $\Bbbk [\Gamma ]$ du groupe $\Gamma$.

Reprenons les notations $L (Q)$ et $D(Q)$ de 5.3, et soit $E$
l'ensemble des semi-invariants non nuls de $B(Q)$. Si $f \in E$ a pour
poids $n_{1} \mu_{1} + \cdots + n_{\ell} \mu_{\ell}$, avec $n_{1},
\dots , n_{\ell} \in \mathbb{Z}$, alors $f\alpha_{1}^{-n_{1}} \cdots
\alpha_{\ell}^{-n_{\ell}} \in L(Q)$ est central dans $L (Q)$. C'est
aussi un \'{e}l\'{e}ment de $B(Q)_{e}$. On voit donc que la
sous-alg\`{e}bre de $B(Q)_{e}$ engendr\'{e}e par $e_{1}, \dots , e_{n},
e^{-1}$ est engendr\'{e}e par $E$ et le centre $C$ de $B(Q)_{e}$. On a
alors $B(Q)_{E} = \big( B(Q)_{e}\big)_{C}$. Compte tenu de 5.4, le
centre de $B(Q)_{E}$ est $D = D(Q)$. On peut donc consid\'{e}rer
$B(Q)_{E}$ comme une P-alg\`{e}bre sur le corps $D$. 

\begin{th*}{Lemme}La $D$-sous-alg\`{e}bre de $B(Q)_{E}$ engendr\'{e}e
par $\Gamma$ s'identifie \`{a} la $D$-alg\`{e}bre $D[\Gamma ]$ du
groupe $\Gamma$. \end{th*}

\begin{Preuve}Soient $\theta_{1}, \dots , \theta_{r} \in D \backslash
\{ 0\}$ et $\gamma_{1}, \dots , \gamma_{r}$ des \'{e}l\'{e}ments
distincts de $\Gamma$ tels que $\theta_{1} \gamma_{1} + \cdots +
\theta_{r} \gamma_{r} = 0$. Alors les $\theta_{i} \gamma_{i}$ sont des
semi-invariants de $B(Q)_{E}$ de poids deux \`{a} deux
distincts. D'apr\`{e}s 5.2, on a $\theta_{i} \gamma_{i} = 0$ pour tout
$i$. Contradiction. \end{Preuve}

\subsection Dans la suite, la $\Bbbk$-sous-alg\`{e}bre
(resp. $D$-sous-alg\`{e}bre) de $B(Q)_{E}$ engendr\'{e}e par $\Gamma$
est encore not\'{e}e $\Bbbk [\Gamma ]$ (resp. $D [\Gamma ]$). On fixe
un suppl\'{e}mentaire $\mathfrak{a}$ de $\mathfrak{n}$ dans
$\mathfrak{g}$. 

\begin{th*}{Lemme}{\em (i)} La dimension $m$ de $\mathfrak{a}$ est
\'{e}gale \`{a} celle du sous-espace de $\mathfrak{g}^{*}$
engendr\'{e} par les \'{e}l\'{e}ments de $G$.

{\em (ii)} Le $D$ sous-espace de $B(Q)_{E}$ engendr\'{e} par
$\mathfrak{a}$ a pour dimension $m$. \end{th*}

\begin{Preuve}(i) Comme on l'a vu en 9.1 :
$$
\mathfrak{n} = \operatornamewithlimits{\textstyle \bigcap}_{\lambda
\in G} \ker \lambda .
$$
D'o\`{u} l'assertion. 

(ii) Soit $\mathscr{S} = (\nu_{1}, \dots , \nu_{m})$ une base du
sous-espace de $\mathfrak{g}^{*}$ engendr\'{e} par $G$, et form\'{e}e
d'\'{e}l\'{e}ments de $G$. Pour $1 \leqslant i \leqslant m$, soit
$a_{i} \in B(Q)_{E}$ un semi-invariant de poids $\nu_{i}$. Notons
$(x_{1}, \dots , x_{m})$ la base duale de $\mathscr{S}$ dans
$\mathfrak{a}$. On a ainsi $\{ x_{i}, a_{j}\} = \delta_{ij} a_{j}$
pour $1 \leqslant i, j \leqslant m$. 

Si $\theta_{1}, \dots , \theta_{m} \in D$ v\'{e}rifient $\theta_{1}
x_{1} + \cdots + \theta_{m} x_{m} = 0$, on obtient
$$
\theta_{i} = \{ \theta_{1} x_{1} + \cdots + \theta_{m} x_{m}, a_{i}\}
= 0
$$
pour $1 \leqslant i \leqslant m$. On a obtenu (ii). \end{Preuve}

\subsection Dans la suite, on note $\mathscr{G}$ (resp. $\mathscr{A}$)
le $D$-sous-espace vectoriel de $B(Q)_{E}$ engendr\'{e} par
$\mathfrak{g}$ (resp. $\mathfrak{a}$).

\subsection Soit $\lambda$ un poids de $\mathfrak{n}$ (pour la
repr\'{e}sentation adjointe de $\mathfrak{g}$ dans
$\mathfrak{n}$). D'apr\`{e}s 5.2, il existe une base $(v_{1,\lambda},
v_{2, \lambda}, \dots , v_{r_{\lambda}, \lambda})$ de
$\mathfrak{n}^{\lambda}$ telle que, pour $x \in \mathfrak{g}$ et $1
\leqslant i \leqslant r_{\lambda}$, on ait 
$$
[x, v_{i,\lambda}] \in \lambda (x) v_{i,\lambda} + \Bbbk v_{i-1,
\lambda} + \cdots + \Bbbk v_{1, \lambda}.
$$
On a ainsi $v_{1, \lambda} \in E \subset D [\Gamma ]$, et il vient 
$$
\{ x, v_{1,\lambda}^{-1}v_{i,\lambda}\} \in \Bbbk v_{1,
\lambda}^{-1}v_{i-1,\lambda} + \cdots + \Bbbk v_{1,\lambda}^{-1}
v_{1, \lambda}.
$$

Notons $A$ la $D$-sous-alg\`{e}bre de $B(Q)_{E}$ engendr\'{e}e par la
r\'{e}union des ensembles $\{ v_{1,\lambda}^{-1}v_{i, \lambda}\, ; \,
1 \leqslant i \leqslant r_{\lambda}\}$ lorsque $\lambda$ parcourt
l'ensemble des poids de $\mathfrak{n}$. Si $x \in \mathfrak{g}$,
$\varepsilon_{Q}(x)$ induit une P-d\'{e}rivation localement nilpotente
de $A$. D'autre part, $A$ est une P-alg\`{e}bre commutative puisque
$\mathfrak{n}$ est ab\'{e}lien. Enfin, $D$ est l'ensemble des $p \in
A$ qui v\'{e}rifient $\varepsilon_{Q}(x) (p) = 0$ pour tout $x \in
\mathfrak{g}$. La commutativit\'{e} de $\mathfrak{n}$ implique aussi
que les restrictions des $\varepsilon_{Q}(x)$ \`{a} $A$ commutent deux
\`{a} deux. 

D'apr\`{e}s 3.8, il existe un $D$-sous-espace vectoriel $W$ de
dimension finie de $A$ v\'{e}rifiant les conditions suivantes :

(i) $A$ est isomorphe \`{a}\ l'alg\`{e}bre sym\'{e}trique $\S (W)$ de
$W$.

(ii) Si $w \in W$ et $y \in \mathscr{G}$, on a $\{ y, w\} \in D$.

(iii) Les seuls id\'{e}aux de $\S (W)$ stables par les $d_{y}$, $y \in
\mathscr{G}$, sont $\{ 0\}$ et $\S (W)$.

\subsection R\'{e}sumons une partie des r\'{e}sultats obtenus
jusqu'ici. 

$\bullet$ La P-alg\`{e}bre $B(Q)_{E}$ est engendr\'{e}e par les
P-sous-alg\`{e}bres $D [\Gamma ]$ et $\S (W)$, ainsi que par le
$D$-sous-espace $\mathscr{A}$.

$\bullet$ Pour tout $y \in \mathscr{A}$, $d_{y}| D[\Gamma ]$ est une
P-d\'{e}rivation semi-simple de $D [\Gamma ]$. De plus, si $y,z \in
\mathscr{A}$, on a $d_{y} {\circ} d_{z} | D [\Gamma ] = d_{z} {\circ}
d_{y}| D [\Gamma ]$.

$\bullet$ Soient $A_{1}$ l'image de $\S (\mathfrak{n})$ dans
$B(Q)_{E}$ et $A$ la P-sous-alg\`{e}bre (sur le corps $D$)
engendr\'{e}e par $(A_{1})_{e}$. Alors $A$ est engendr\'{e}e par $D
[\Gamma ]$ et $\S (W)$. En outre, si $y \in \mathscr{A}$, $d_{y}| \S
(W)$ est une d\'{e}rivation localement nilpotente de $\S
(W)$. D'autre part, si $y,z \in \mathscr{A}$, on a $\{ y,z\} \in
A$. 

\subsection Notons $R$ la P-alg\`{e}bre $\S (W) \otimes_{D} D
[\Gamma ]$ (voir 2.8) sur le corps $D$. On fait op\'{e}rer
$\mathscr{A}$ sur $R$ en convenant que, si $y \in \mathscr{A}$, $a
\in \S (W)$ et $g \in D [\Gamma ]$, on a :
$$
y.(a\otimes g) = \{ y,a\} \otimes g + a \otimes \{ y,g\} .
$$

\begin{th*}{\bf Lemme}Les seuls {\em P}-id\'{e}aux de $R$ qui sont
stables sous l'action de $\mathscr{A}$ sont $\{ 0\}$ et $R$. \end{th*}

\begin{Preuve}On remarque que l'op\'{e}ration de $\mathscr{A}$ sur $R$
est localement finie. 

Soit $g \in \Gamma$. Il existe une forme $D$-lin\'{e}aire $\mu$ sur
$\mathscr{A}$ telle que $\{ y,g\} = \mu (y) g$ pour tout $y \in
\mathscr{A}$. En outre, \`{a} des \'{e}l\'{e}ments distincts de
$\Gamma$ correspondent des formes lin\'{e}aires distinctes (voir 9.1). 

Avec ces notations, l'espace $R^{\mu}$ est \'{e}gal \`{a} $\S (W)
\otimes g$ car, $\mathscr{A}$ op\`{e}re de mani\`{e}re localement
nilpotente sur $\S (W)$. 

Soit $I$ un P-id\'{e}al de $R$ stable pour l'action de
$\mathscr{A}$. D'apr\`{e}s ce qui pr\'{e}c\`{e}de, on a :
$$ 
I = \operatornamewithlimits{\textstyle \sum}_{g \in \Gamma} I \cap
\big( \S (W) \otimes g\big) .
$$
Supposons $I$ non nul. Il existe $g \in \Gamma$ tel que $I \cap \big(
\S (W) \otimes g\big) \ne \{ 0\}$ et, $I$ \'{e}tant un id\'{e}al de
$R$, il vient $I \cap \S (W) \ne \{ 0\}$. Le seul P-id\'{e}al non nul
de $\S (W)$ stable sous l'action de $\mathscr{A}$ est $\S (W)$
d'apr\`{e}s 3.8. On a donc obtenu le r\'{e}sultat. \end{Preuve}

\begin{th}{\bf Lemme}La {\em P}-sous-alg\`{e}bre (sur le corps $D$) de
$B(Q)_{E}$ engendr\'{e}e par $\S (W)$ et $D [\Gamma ]$ est isomorphe
\`{a} $\S (W) \otimes_{D} D [\Gamma ]$. \end{th}

\begin{Preuve}Notons $A$ cette P-alg\`{e}bre. Elle est commutative, et
  l'application 
$$
\S (W) \otimes_{D} D[\Gamma ] \to A \ , \ a \otimes g \to ag
$$
est un homomorphisme surjectif de P-alg\`{e}bres. Le noyau de cet
homomorphisme est un P-id\'{e}al stable sous l'action de
$\mathscr{A}$. Le lemme est donc une cons\'{e}quence de
9.6. \end{Preuve}

\subsection Soit $V$ le $D$-espace vectoriel $W \oplus \mathscr{A}$ de
$B(Q)_{E}$. D'apr\`{e}s 9.4, si $v_{1}, v_{2} \in V$, on a $\{
v_{1}, v_{2}\} \in D$, et l'application
$$
\omega \colon V \times V \to \Bbbk \ , \ (v_{1}, v_{2}) \to \{ v_{1}, v_{2}\}
$$
est une forme bilin\'{e}aire altern\'{e}e sur $V$. D'apr\`{e}s 9.6, le
seul \'{e}l\'{e}ment $v$ de $W$ v\'{e}rifiant $\{ v , \mathscr{A}\} =
\{ 0\}$ est l'\'{e}l\'{e}ment nul. Par suite, le noyau $V^{\omega}$ de
$\omega$ v\'{e}rifie $V^{\omega} \subset \mathscr{A}$.

Pour $v \in V$ et $g \in \Gamma$, on a 
$$
\{ g, v\} = \lambda_{g}(v) g,
$$
o\`{u} $\lambda_{g}$ est une forme lin\'{e}aire sur $V^{*}$. On a
$$ 
W = \operatornamewithlimits{\textstyle \bigcap}_{g \in \Gamma} \ker
\lambda_{g}
$$
et, si $g_{1}, g_{2}$ sont des \'{e}l\'{e}ments distincts de $\Gamma$,
alors $\lambda_{g_{1}} \ne \lambda_{g_{2}}$. On en d\'{e}duit que $G =
\{ \lambda_{g} \, ; \, g \in \Gamma\}$ est un sous-groupe de $V^{*}$
isomorphe \`{a} $\Gamma$.

On a vu que $V^{G} = W$ et, comme $V^{\omega} \subset \mathscr{A}$, il
vient $V^{G} \cap V^{\omega} = \{ 0\}$. Ainsi, la P-alg\`{e}bre
$\mathscr{B}_{D}(V, \omega , G)$ construite comme en 7.1 est simple
d'apr\`{e}s 7.5.

Comme $B(Q)_{E}$ est engendr\'{e}e par $V$ et par $D [\Gamma ]$, et
que $\mathscr{B}_{D}(V, \omega , G)$ est simple, il r\'{e}sulte de la
proposition 7.3 que les P-alg\`{e}bres (sur le corps $D$) $B(Q)_{E}$
et $\mathscr{B}_{D}(V,\omega , G)$ sont isomorphes.

\section{Structure des quotients premiers}

\chead[\sl Alg\`{e}bres de Poisson et alg\`{e}bres de Lie
r\'{e}solubles]{\sl 10~~Structure des quotients premiers}

\subsection Dans tout ce paragraphe, $\mathfrak{g}$ est une
alg\`{e}bre de Lie r\'{e}soluble. 

Soit $Q$ un id\'{e}al premier $\mathfrak{g}$-stable de $\S
(\mathfrak{g})$. On va s'int\'{e}resser \`{a} l'alg\`{e}bre de Poisson
$\S (\mathfrak{g}) / Q$ ; on peut donc supposer que $\mathfrak{g} \cap
Q = \{ 0\}$, et on note $\pi$ la surjection canonique de $\S
(\mathfrak{g})$ dans $B(Q, \mathfrak{g}) = \S (\mathfrak{g})/ Q$. 

Notons $\mathfrak{n}$ le plus grand id\'{e}al nilpotent de
$\mathfrak{g}$, et soit $\mathfrak{a}$ un suppl\'{e}mentaire de
$\mathfrak{n}$ dans $\mathfrak{g}$. Si $R = \S (\mathfrak{n}) \cap Q$,
alors $\pi (\big(\S (\mathfrak{n})\big)$ s'identifie \`{a}
$B(R,\mathfrak{n})$.

Supposons $B(R, \mathfrak{n}) \ne Y(R, \mathfrak{n})$. D'apr\`{e}s
8.2, il existe $e \in B(R, \mathfrak{n})$ et $n \in \mathbb{N}^{*}$
v\'{e}rifiant les conditions suivantes :

(i) $e$ est vecteur propre pour l'action de $\mathfrak{g}$ dans $B(R,
\mathfrak{n})$. 

(ii) L'alg\`{e}bre $Y(R, \mathfrak{n})_{e}$ est de type fini.

(iii) Les P-alg\`{e}bres $B(R, \mathfrak{n})_{e}$ et $Y(R,
\mathfrak{n})_{e} \otimes_{\Bbbk} B_{n}$ sont isomorphes.

La P-alg\`{e}bre $B(Q, \mathfrak{g})_{e}$ est ainsi engendr\'{e}e par
$Y (R, \mathfrak{n})_{e} \otimes_{\Bbbk} B_{n}$ et $\pi
(\mathfrak{a})$. 

\subsection D'apr\`{e}s 3.4,\,(iv), il existe une application
lin\'{e}aire $\rho \colon \mathfrak{a} \to B(R, \mathfrak{n})_{e}$
telle que, pour tout $x \in \mathfrak{a}$, $\pi (x) - \rho (x)$
commute \`{a} $B_{n} \subset B(R,\mathfrak{n})_{e}$. Notons
$\mathscr{A}$ le sous-espace de $B(Q, \mathfrak{g})_{e}$ engendr\'{e}
par les $\pi (x) - \rho (x)$ pour $x \in \mathfrak{a}$. Comme
$[\mathfrak{g}, \mathfrak{g}] \subset \mathfrak{n}$, on a
$\{\mathscr{A}, \mathscr{A}\} \subset B(R, \mathfrak{n})_{e}$, donc
$\{\mathscr{A}, \mathscr{A}\} \subset Y(Q, \mathfrak{n})_{e}$ puisque
$\{ \mathscr{A}, \mathscr{A}\}$ commute \`{a} $B_{n}$. 

Soit $C$ la P-sous-alg\`{e}bre de $B(Q,\mathfrak{g})_{e}$
engendr\'{e}e par $Y(R, \mathfrak{n})_{e}$ et $\mathscr{A}$. Alors
$B(Q,\mathfrak{g})_{e}$ est engendr\'{e} par $B_{n}$ et
$C$. D'apr\`{e}s 3.5,\,(i), les P-alg\`{e}bres $B(Q,\mathfrak{g})_{e}$
et $C \otimes_{\Bbbk} B_{n}$ sont isomorphes. 

L'alg\`{e}bre de Lie $\mathfrak{g}$ \'{e}tant r\'{e}soluble, $Y(R,
\mathfrak{n})_{e}$ est une r\'{e}union croissante de
$\mathfrak{g}$-modules de dimension finie (pour la repr\'{e}sentation
induite par la repr\'{e}sentation adjointe de $\mathfrak{g}$). D'autre
part, si $x \in \mathfrak{a}$ et $u \in Y(R, \mathfrak{n})_{e}$, on a
:
$$
\{ \pi (x) - \rho (x), u\} = \{ \pi (x), u\} .
$$
Comme on a d\'{e}j\`{a} vu que $\{ \mathscr{A}, \mathscr{A}\} \subset
Y(R, \mathfrak{n})_{e}$ et que $Y(R, \mathfrak{n})_{e}$ est une
alg\`{e}bre de type fini, on voit qu'il existe un sous-espace $V$ de
dimension dinie de $Y(R, \mathfrak{n})_{e}$ v\'{e}rifiant les
conditions suivantes :

(i) $\{ \mathscr{A}, \mathscr{A}\} \subset V$.

(ii) $V$ engendre l'alg\`{e}bre $Y(R, \mathfrak{n})_{e}$. 

Posons $\mathfrak{h} = \mathscr{A} + V$. Il est imm\'{e}diat que l'on
d\'{e}finit une structure d'alg\`{e}bre de Lie r\'{e}soluble sur
$\mathfrak{h}$ en posant $[u,v] = \{ u, v\}$ pour $u,v \in
\mathfrak{h}$. Enfin, il est clair que la P-alg\`{e}bre $C$
pr\'{e}c\'{e}dente est le quotient de $\S (\mathfrak{h})$ par un
id\'{e}al premier $\mathfrak{h}$-stable $T$ de $\S
(\mathfrak{h})$. Les P-alg\`{e}bres $B(Q,\mathfrak{g})_{e}$ et $B(T,
\mathfrak{h}) \otimes_{\Bbbk} B_{n}$ sont alors isomorphes et, comme
en 10.1, on peut supposer que $\mathfrak{h} \cap T = \{ 0\}$. 

\subsection Compte-tenu des alin\'{e}as pr\'{e}c\'{e}dents, on voit
que l'on se ram\`{e}ne au cas suivant. Il existe un semi-invariant
non nul $e$ de $B(Q, \mathfrak{g})$, une alg\`{e}bre de Lie
r\'{e}soluble $\mathfrak{h}$, un id\'{e}al premier
$\mathfrak{h}$-stable de $\S (\mathfrak{h})$ et un entier positif ou
nul $n$ v\'{e}rifiant les conditions suivantes :

(i) Les P-alg\`{e}bres $B(Q, \mathfrak{g})_{e}$ et $B(T, \mathfrak{h})
\otimes_{\Bbbk} B_{n}$ sont isomorphes.

(ii) Le plus grand id\'{e}al nilpotent de $\mathfrak{h}$ est
commutatif. 

\begin{th}{\bf Th\'{e}or\`{e}me}On suppose que le corps $\Bbbk$ est
alg\'{e}briquement clos. Soient $\mathfrak{g}$ une alg\`{e}bre de
Lie r\'{e}soluble et $Q$ un id\'{e}al pemier $\mathfrak{g}$-stable
de $\S (\mathfrak{g})$. On note $B(Q)$ la {\em P}-alg\`{e}bre $\S
(\mathfrak{g}) / Q$, $E$ l'ensemble des semi-invariants non nuls de
$B(Q)$, et $D$ le corps des invariants de $\Fract B(Q)$. Sur le
corps $D$, la {\em P}-alg\`{e}bre $B(Q)_{E}$ est isomorphe \`{a} une
{\em P}-alg\`{e}bre simple $\mathscr{B}_{D}(V, \omega , G)$
d\'{e}finie comme en {\em 7.2}.\end{th}

\begin{Preuve}C'est une cons\'{e}quence de 9.8, 10.3, et
3.4,\,(iii). \end{Preuve} 

\subsection On d\'{e}signe par $SD(Q, \mathfrak{g})$ ou $SD (Q)$ le
corps des fractions de $Y(Q)$. On note $\mathscr{V}(Q)$ la
vari\'{e}t\'{e} des z\'{e}ros de $Q$ dans $\mathfrak{g}^{*}$ et, si
$H$ est un sous-groupe alg\'{e}brique du groupe des automorphismes de
$\mathfrak{g}$, $m_{H}(Q)$ est la dimension maximale des $H$-orbites
dans $\mathscr{V}(Q)$.

Avec les hypoth\`{e}ses et notations de 10.4, on peut interpr\'{e}ter
certains entiers li\'{e}s \`{a} la P-alg\`{e}bre $\mathscr{B}_{D}(V,
\omega , G)$ (voir 7.8) en fonction de $Q$. Les r\'{e}sultats qui
suivent sont prouv\'{e}s dans le paragraphe 3 de la partie III de [9].

\begin{th*}{\bf Proposition}Soient $\mathfrak{g}, Q$ et
$\mathscr{B}_{D}(V, \omega , G)$ comme en {\em 10.4}, $\Gamma$ le
groupe al\-g\'{e}\-bri\-que adjoint de $\mathfrak{g}$, et $\Gamma_{0}$
l'intersection des noyaux des caract\`{e}res rationnels de
$\Gamma$. Alors :

{\em (i)} $\deg \tr_{\Bbbk} D(Q) = \dim_{\Bbbk} \mathscr{V}(Q) -
m_{\Gamma}(Q)$.

{\em (ii)} $\deg \tr_{\Bbbk} SD(Q) = \dim_{\Bbbk} \mathscr{V} (Q) -
m_{\Gamma_{0}} (Q)$.

{\em (iii)} $\rg (G) = \deg \tr_{\Bbbk} SD (Q) - \deg \tr_{\Bbbk} D(Q)$.

{\em (iv)} $\dim_{D(Q)} V = \dim_{\Bbbk} \mathscr{V}(Q) - \deg
\tr_{\Bbbk} SD (Q)$.

\end{th*}

\section{Cas alg\'{e}brique}

\chead[\sl Alg\`{e}bres de Poisson et alg\`{e}bres de Lie
r\'{e}solubles]{\sl 10~~Cas alg\'{e}brique}

\subsection On d\'{e}signe toujours par $Q$ un id\'{e}al premier
$\mathfrak{g}$-stable de $\S (\mathfrak{g})$, et on conserve les
notations du paragraphe 5.

Si $\mathscr{J}(\mathfrak{g}) = \{ \nu_{1}, \dots , \nu_{m}\}$, il
r\'{e}sulte de 5.4,\,(ii) et 6.2 que
$$
\Lambda'(Q) \subset \mathbb{Z} \nu_{1} + \cdots + \mathbb{Z}
\nu_{m}. 
$$
Le corps $\Bbbk$ \'{e}tant de caract\'{e}ristique nulle, $\Lambda'(Q)$
est un $\mathbb{Z}$-module libre ; soit $(\mu_{1}, \dots ,
\mu_{\ell})$ une base de ce module libre.

A nouveau d'apr\`{e}s 5.4,\,(ii), pour $1 \leqslant i \leqslant \ell$,
on peut \'{e}crire
$$
\mu_{i} = \operatornamewithlimits{\textstyle \sum}_{j=1}^{s} n_{ij}
\lambda_{j}, 
$$
o\`{u} $\lambda_{1}, \dots , \lambda_{s} \in \Lambda (Q)$, et o\`{u}
les $n_{ij}$ sont des entiers.

Pour $1 \leqslant j \leqslant s$, soit $e_{j} \in E (Q ,
\mathfrak{g})$ un \'{e}l\'{e}ment non nul de poids
$\lambda_{j}$. Posons $e = e_{1} \cdots e_{s}$. L'id\'{e}al $Q$
\'{e}tant premier, on a $e^{n} \notin Q$ pour tout $n \in \mathbb{N}$,
et on peut former la P-alg\`{e}bre $B(Q, \mathfrak{g})_{e}$. 

On note $Y(B_{e})$ le centre de la P-alg\`{e}bre
$B(Q,\mathfrak{g})_{e}$ et, avec les notations de 5.2, on pose :
$$
SY (B_{e}) = \operatornamewithlimits{\textstyle \sum}_{\lambda \in
\mathfrak{g}^{*}} \big( B((Q,\mathfrak{g})_{e}\big)_{\lambda} .
$$

Si $1 \leqslant i \leqslant \ell$, on fixe un mon\^{o}me $a_{i}$ en
les $e_{j}, e_{j}^{-1}$, de poids $\mu_{i}$, et on note $A$ le
sous-groupe multiplicatif de $B(Q,\mathfrak{g})_{e}$ engendr\'{e} par
$a_{1}, \dots , a_{\ell}$.

\begin{th*}{\bf Lemme}{\em (i)} Soient $u = a_{1}^{m_{1}} \cdots
a_{\ell}^{m_{\ell}}$ et $v = a_{1}^{n_{1}} \cdots
a_{\ell}^{m_{\ell}}$, o\`{u} $m_{1}, n_{1}, \dots , m_{\ell} ,
n_{\ell}$ sont des entiers. Si $u$ et $v$ ont m\^{e}me poids, alors
$m_{i} = n_{i}$ pour $1 \leqslant i \leqslant \ell$.

{\em (ii)} Si $a \in A$ a pour poids $\lambda \in \mathfrak{g}^{*}$,
l'ensemble des \'{e}l\'{e}ments semi-invariants de poids $\lambda$ de
$B(Q, g)_{e}$ est \'{e}gal \`{a} $a Y(B_{e})$.

{\em (iii)} Les \'{e}l\'{e}ments $a_{1}, \dots , a_{\ell}$ sont
alg\'{e}briquement ind\'{e}pendants sur $Y(B_{e})$.

{\em (iv)} On a $SY (B_{e}) = Y (B_{e}) \{ a_{1}, a_{1}^{-1}, \dots,
a_{\ell}, a_{\ell}^{-1}\}$. \end{th*}

\begin{Preuve}(i) C'est imm\'{e}diat puisque $(\mu_{1}, \dots ,
\mu_{\ell})$ est une base du $\mathbb{Z}$-module libre $\Lambda' (Q)$.

(ii) Si $b \in B(Q,\mathfrak{g})_{e}$ a m\^{e}me poids que
$a$, alors $a^{-1}b \in Y(B_{e})$. D'o\`{u} l'assertion.

(iii) Si $i_{1}, \dots , i_{\ell} \in \mathbb{N}$ et $\alpha \in
Y(B_{e})$, alors $\alpha a_{1}^{i_{1}} \cdots a_{\ell}^{i_{\ell}}$ est
un semi-invariant de poids $i_{1} \mu_{1} + \cdots + i_{\ell}
\mu_{\ell}$. Par suite, si 
$$
P = \operatornamewithlimits{\textstyle \sum} \alpha_{i_{1}, \dots ,
i_{\ell}} X_{1}^{i_{1}} \cdots X_{\ell}^{i_{\ell}} \in
\Bbbk [X_{1}, \dots , X_{\ell}]
$$
v\'{e}rifie $P(a_{1}, \dots , a_{\ell}) = 0$, il r\'{e}sulte de 5.2 et
du fait que $(\mu_{1}, \dots , \mu_{\ell})$ est une base de
$\Lambda'(Q)$, que $P = 0$. 

(iv) En tant qu'alg\`{e}bres associatives, $SY (B_{e}) = Y(B_{e})
[a_{1}, a_{1}^{-1}, \dots , a_{\ell}, a_{\ell}^{-1}]$ d'apr\`{e}s ce
qui pr\'{e}c\`{e}de. Le r\'{e}sultat est donc cons\'{e}quence de
5.10. \end{Preuve}

\subsection Soit $\mathfrak{h}$ une alg\`{e}bre de Lie. On dit que
$\mathfrak{h}$ est alg\'{e}brique si elle est isomorphe \`{a}
l'alg\`{e}bre de Lie d'une groupe alg\'{e}brique. Pour ce qui concerne
les alg\`{e}bres de Lie alg\'{e}briques, le lecteur pourra se reporter
\`{a} [10]. Nous allons seulement rappeler ici quelques r\'{e}sultats
sur les alg\`{e}bres de Lie r\'{e}solubles alg\'{e}briques ; ils
proviennent de [4] et [8].

\begin{th}{\bf Proposition}Soit $\mathfrak{g}$ une alg\`{e}bre de Lie
r\'{e}soluble alg\'{e}brique.

{\em (i)} Toute image de $\mathfrak{g}$ par homomorphisme est
alg\'{e}brique. 

{\em (ii)} On a $\mathfrak{g} = \mathfrak{s} \oplus \mathfrak{n}$,
o\`{u} $\mathfrak{n}$ est la plus grand id\'{e}al nilpotent de
$\mathfrak{g}$, et o\`{u} $\mathfrak{s}$ est une sous-alg\`{e}bre de
Lie ab\'{e}lienne de $\mathfrak{g}$ qui op\`{e}re de mani\`{e}re
semi-simple sur $\mathfrak{n}$.

{\em (iii)} Soit $\mathscr{J}(\mathfrak{g})$ comme en {\em 6.1}. Le
rang du $\mathbb{Z}$-sous-module de $\mathfrak{g}^{*}$ engendr\'{e}
par $\mathscr{J}(\mathfrak{g})$ est \'{e}gal \`{a} la dimension du
$\Bbbk$-sous-espace vectoriel de $\mathfrak{g}^{*}$ engendr\'{e} par
$\mathscr{J}(\mathfrak{g})$. \end{th}

\subsection Avec les hypoth\`{e}ses de 11.3, soit $H$ (resp. $V$) un
$\mathbb{Z}$-sous-module (resp. un $\Bbbk$-sous-espace vectoriel) du
$\mathbb{Z}$-module (resp. de l'espace vectoriel) engendr\'{e} par
$\mathscr{J}(\mathfrak{g})$. Si l'espace vectoriel engendr\'{e} par
$H$ est $V$, alors le rang de $H$ est \'{e}gal \`{a} la dimension de
$V$.

\subsection On suppose dor\'{e}navant que $\mathfrak{g}$ est une
alg\`{e}bre de Lie r\'{e}soluble al\-g\'{e}\-bri\-que, et on \'{e}crit
$\mathfrak{g} = \mathfrak{s} \oplus \mathfrak{n}$ comme en 11.3. On
conserve les notations de 11.1 et du paragraphe 5, et on note
$\widehat{\mathfrak{g}}$ pour $\widehat{\mathfrak{g}}_{Q}$.

Pour \'{e}tudier les localis\'{e}s de la P-alg\`{e}bre $B(Q,
\mathfrak{g})$ on peut supposer, d'apr\`{e}s 11.3,\,(i) que $Q \cap
\mathfrak{g} = \{ 0\}$, ce que nous ferons d\'{e}sormais. On
d\'{e}signe par $\{ \mu_{1}, \dots ,\mu_{\ell}\}$ une base du
$\mathbb{Z}$-module $\Lambda'(Q)$ ; d'apr\`{e}s 11.3, c'est aussi une
base du $\Bbbk$-espace vectoriel engendr\'{e} par $\Lambda'(Q)$. On a
:
$$
\widehat{\mathfrak{g}} = \operatornamewithlimits{\textstyle
\bigcap}_{i=1}^{\ell} \ker \mu_{i}.
$$
Ecrivons $\mathfrak{g} = \mathfrak{s} \oplus \mathfrak{n}$ (notations
de 11.3,\,(ii)). On a $\mathfrak{n} \subset
\widehat{\mathfrak{g}}$. On d\'{e}signe par $\mathfrak{t}$ une
sous-alg\`{e}bre de Lie ab\'{e}lienne de $\mathfrak{s}$ qui
v\'{e}rifie $\mathfrak{g} = \mathfrak{t} \oplus
\widehat{\mathfrak{g}}$, et qui op\`{e}re de mani\`{e}re semi-simple
sur $\widehat{\mathfrak{g}}$. Pour tout poids $\lambda$ d'un
semi-invariant de $B(Q,\mathfrak{g})$, on a $\widehat{\mathfrak{g}}
\subset \ker \lambda$. On en d\'{e}duit facilement que 
$\{ \mu_{1}|\mathfrak{t}, \dots , \mu_{\ell}|\mathfrak{t}\}$ est une
base de $\mathfrak{t}^{*}$. 

\begin{th}{\bf Lemme}On a $SY (Q, \mathfrak{g}) = Y (\widehat{Q},
\widehat{\mathfrak{g}}) = SY (\widehat{Q},
\widehat{\mathfrak{g}})$. \end{th}

\begin{Preuve}On a vu en 5.13 que $SY (Q, \mathfrak{g}) \subset Y
(\widehat{Q}, \widehat{\mathfrak{g}}) = SY (\widehat{Q},
\widehat{\mathfrak{g}})$. D'apr\`{e}s les hypoth\`{e}ses,
$\mathfrak{t}$ op\`{e}re de mani\`{e}re semi-simple sur
$Y(\widehat{Q}, \widehat{\mathfrak{g}})$. Comme $\{
\widehat{\mathfrak{g}} , Y (\widehat{Q}, \widehat{\mathfrak{g}})\} =
\{ 0\}$, on voit que tout \'{e}l\'{e}ment de $Y(\widehat{Q},
\widehat{\mathfrak{g}})$ est somme de vecteurs propres pour l'action
de $\mathfrak{g}$. D'o\`{u} $Y (\widehat{Q}, \widehat{\mathfrak{g}})
\subset SY (Q, \mathfrak{g})$. \end{Preuve}

\subsection Il r\'{e}sulte de 8.2 qu'il existe $e \in Y
(\widehat{Q}, \widehat{\mathfrak{g}}) \backslash \{ 0\}$ et $p \in
\mathbb{N}$ v\'{e}rifiant les conditions suivantes :

(i) L'alg\`{e}bre $Y (\widehat{Q}, \widehat{\mathfrak{g}})_{e}$ est de
type fini et $e$ est un semi-invariant de $B(Q, \mathfrak{g})$.

(ii) Les P-alg\`{e}bres $B(\widehat{Q}, \widehat{\mathfrak{g}})_{e}$
et $Y(\widehat{Q}, \widehat{\mathfrak{g}})_{e} \otimes_{\Bbbk} B_{p}$
sont isomorphes, et on peut supposer que des g\'{e}n\'{e}rateurs $x_{1},
y_{1}, \dots , x_{p}, y_{p}$ de $B_{p} \subset B(\widehat{Q},
\widehat{\mathfrak{g}})_{e}$ sont des vecteurs propres pour l'action de
$\mathfrak{t}$ dans $B (\widehat{Q}, \widehat{\mathfrak{g}})_{e}$. 

Pour $1 \leqslant i \leqslant p$, notons $\nu_{i}, \theta_{i}$ des
formes lin\'{e}aires sur $\mathfrak{t}$ telles que
$$
\{ t, x_{i}\} = \nu_{i}(t) x_{i} \ , \ \{ t, y_{i}\} = \theta_{i}(t) y_{i}
$$
pour tout $t \in \mathfrak{t}$. Comme $\{ x_{i}, y_{i}\} = 1$, on
obtient facilement $\theta_{i} = - \nu_{i}$. Si $t \in \mathfrak{t}$,
posons alors 
$$ 
t' = t + \nu_{1}(t) x_{1} y_{1} + \cdots + \nu_{p}(t) x_{p}y_{p}.
$$ 
Un calcul imm\'{e}diat montre que, pour $s,t \in \mathfrak{t}$ et
$1 \leqslant i \leqslant p$, on a :
$$
\{ s', t'\} = \{ t', x_{i}\} = \{ t', y_{i}\} = 0.
$$

La P-alg\`{e}bre $B(Q,\mathfrak{g})_{e}$ est engendr\'{e}e par
$B(\widehat{Q}, \widehat{\mathfrak{g}})_{e}$ et les $t \in
\mathfrak{t}$, c'est-\`{a}-dire aussi par $B(\widehat{Q},
\widehat{\mathfrak{g}})_{e}$ et les $t'$, avec $t \in
\mathfrak{t}$. Compte tenu de ce qui pr\'{e}c\`{e}de et de
5.13,\,(ii), si $(z_{1}, \dots , z_{\ell})$ est une base de
$\mathfrak{t}$, la P-alg\`{e}bre $B(Q,\mathfrak{g})_{e}$ est isomorphe
\`{a} la P-alg\`{e}bre
$$
\big( \cdots \big( \big( Y(\widehat{Q},
\widehat{\mathfrak{g}})_{e}\big)_{\delta_{1}} \{ X_{1}\} \big) \cdots
\big)_{\delta_{r}} \{ X_{r}\} \otimes_{\Bbbk} B_{p},
$$
o\`{u} $\delta_{i}$ est la P-d\'{e}rivation induite par $z_{i}$.

D'apr\`{e}s 11.1, il existe $f \in E(Q)$ tel que l'on ait 
$$
SY (B_{f}) = Y (B_{f}) \otimes_{\Bbbk} \Bbbk [a_{1}, a_{1}^{-1}, \dots
, a_{\ell}, a_{\ell}^{-1}],
$$
o\`{u} $a_{i}$ est un semi-invariant de poids $\mu_{i}$. Quitte \`{a}
remplacer $f$ par $ef$, on peut alors supposer que l'on a
$$
SY (B_{e}) = Y(B_{e}) \otimes_{\Bbbk} \Bbbk [a_{1}, a_{1}^{-1}, \dots
, a_{\ell}, a_{\ell}^{-1}] .
$$

On a d\'{e}j\`{a} remarqu\'{e} que $(\mu_{1}| \mathfrak{t}, \dots ,
\mu_{\ell} | \mathfrak{t})$ est une base de $\mathfrak{t}^{*}$. Notons
$(b_{1}, \dots , b_{\ell})$ la base duale dans $\mathfrak{t}$, et
posons $t_{i} = a_{i}^{-1}b_{i}$ pour $1 \leqslant i \leqslant
\ell$. Comme $a_{j}$ est un semi-invariant de poids $\mu_{i}$, on
obtient facilement, pour $1 \leqslant i,j \leqslant \ell$, 
$$
\{ t_{i}, a_{j}\} = \delta_{ij} \ , \ \{ t_{i}, t_{j}\} = 0,
$$
et on a $\{ a_{i}, a_{j}\} = 0$ puisque $a_{i}$ et $a_{j}$ sont
semi-invariants (5.4,\,(ii) et 5.10). On a en particulier
$$
\{ t_{i}, a_{1}^{m_{1}}\cdots a_{\ell}^{m_{\ell}}\} = m_{i}
a_{1}^{m_{1}} \cdots a_{i-1}^{m_{i-1}} a_{i}^{m_{i}-1}
a_{i+1}^{m_{i+1}} \cdots a_{\ell}^{m_{\ell}}
$$
pour $m_{1}, \dots , m_{\ell} \in \mathbb{N}$. 

Comme $SY (B_{e}) = \big( SY(Q)\big)_{e}$ et que $Y(B_{e}) = \big(
Y(Q)\big)_{e}$, ce qui pr\'{e}c\`{e}de prouve que les P-alg\`{e}bres
$\big( \cdots \big( \big( Y(\widehat{Q},
\widehat{\mathfrak{g}})_{e}\big)_{\delta_{1}} \{ X_{1}\} \big) \cdots
\big)_{\delta_{\ell}} \{ X_{\ell}\}$ et $Y(B_{e}) \otimes_{\Bbbk}
B'_{\ell}$ sont isomorphes (o\`{u} $B'_{\ell}$ a la m\^{e}me
signification qu'en 3.10).

On a donc obtenu le r\'{e}sultat suivant :

\begin{th*}{\bf Th\'{e}or\`{e}me}On suppose que $\mathfrak{g}$ est
alg\'{e}brique et que $Q$ est un id\'{e}al premier
$\mathfrak{g}$-stable de $\S (\mathfrak{g})$. Il existe un
semi-invariant non nul $e$ de $B(Q, \mathfrak{g})$ et des entiers $p$
et $\ell$ tels que les {\em P}-alg\`{e}bres $B(Q,\mathfrak{g})_{e}$ et
$Y\big( B(Q, \mathfrak{g})_{e}\big) \otimes_{\Bbbk} B_{p}
\otimes_{\Bbbk} B'_{\ell}$ soient isomorphes. \end{th*}

\begin{th}{\bf Corollaire}Avec les hypoth\`{e}ses pr\'{e}c\'{e}dentes, il
existe un entier $n$ tels que, sur le corps $D(Q, \mathfrak{g})$, les
{\em P}-alg\`{e}bres $\Fract B(Q,\mathfrak{g})$ et $\Fract B_{n}\big(
D(Q, \mathfrak{g})\big)$ soient isomorphes. \end{th}

\begin{th*}{\bf Remarque}{\em Le r\'{e}sultat de 11.8 est l'analogue,
pour les P-alg\`{e}bres, de l'assertion (i) de [4], corollary 6.4.}
\end{th*}

\medskip

\begin{center}
{\titre Bibliographie}
\end{center}

\medskip

[1] {\sc Borho W., Gabriel P., Rentschler R.}, \textit{Primideale in
  Einh\"{u}llenden aufl\"{o}bare Lie Algebren}, Lectures Notes in
  mathematics, 357, Springer Verlag, 1973.

[2] {\sc Dixmier J.}, \textit{Enveloping algebras}, Graduate Studies
en Math., 11, AMS, 1996.

[3] {\sc Krause G. R., Lenagan T. H.}, \textit{Growth of Algebras and
Gelfand-Kirillov Dimension}, Graduate Studies in Math., 22, AMS, 2000.

[4] {\sc Mc\,Connell J.C}, Representations of solvable Lie algebras
and the Gel\-fand-Kirillov conjecture, \textit{Proc. London Math. Soc.},
series 3, 29, 1974, p. 453-484.

[5] {\sc Mc\,Connell J.C}, Representations of solvable Lie algebras II ;
Twisted group-rings, \textit{Ann. Scient. Ec. Norm. Sup.}, s\'{e}rie
4, t.8, 1975, p. 157-178.

[6] {\sc Mc\,Connell J.C}, Representations of solvable Lie algebras
IV : An elementary proof of the $(U/P)_{E}$-structure, \textit{Proc. of
the American Math. Soc.}, Vol. 66, 1, 1977, p. 8.12.

[7] {\sc Mc\,Connell J.C}, Amalgams of Weyl algebras and the
$\mathscr{A}(V, \delta , \Gamma )$ conjecture, \textit{Invent. math.},
92, 1988, p. 163-171.

[8] {\sc Mc\,Connell J.C, Robson J.C}, \textit{Noncommutative
noetherian rings}, Wiley series in pure and applied mathematics,
John Wiley \& Sons, 1987.

[9] {\sc Tauvel P.}, Sur les quotients premiers de l'alg\`{e}bre
enveloppante d'une alg\`{e}bre de Lie r\'{e}soluble,
\textit{Bull. Soc. Math. France}, 10, 1978, p. 177-205.

[10] {\sc Tauvel P., Yu R. W. T.}, \textit{Lie algebras and Algebraic
  Groups}, Springer Monographs in Math., Springer, 2005.

[11] {\sc Vergne M.}, La structure de Poisson sur l'alg\`{e}bre
sym\'{e}trique d'une alg\`{e}bre de Lie nilpotente,
\textit{Bull. Soc. Math. France}, 100, 1984, p. 301-335.

\bigskip

\noindent{\sc UMR 6086 CNRS, D\'{e}partement de Math\'{e}matiques\\
T\'{e}l\'{e}port 2, BP 30179 \\
 Boulevard Marie et Pierre Curie\\
86962 Futuroscope Chasseneuil Cedex \\
France}

\medskip

\textit{Adresses e-mail~:}

tauvel@math.univ-poitiers.fr \ (P. Tauvel)

yuyu@math.univ-poitiers.fr \ (R.W.T Yu)

\end{document}